\documentclass[12pt]{amsart}
\usepackage[margin=1in]{geometry}

\usepackage{amssymb}
\usepackage{hyperref} 
\usepackage{enumitem}
\usepackage{comment}
\usepackage{amsmath}
\usepackage{latexsym}
\usepackage{extarrows}
\usepackage{enumerate}
\usepackage{txfonts}
\usepackage{mathtools}
\usepackage{bm}
\usepackage{tikz-cd}
\usepackage{xcolor}
\IfFileExists{ulem.sty}{\usepackage[normalem]{ulem}}{}

\makeatletter
\@namedef{subjclassname@2020}{%
  \textup{2020} Mathematics Subject Classification}
\makeatother

\usepackage[T1]{fontenc}

\theoremstyle{plain}

\newtheorem{theorem}{Theorem}[section]
\newtheorem{proposition}[theorem]{Proposition}

\newtheorem{lemma}[theorem]{Lemma}

\newtheorem{conjecture}[theorem]{Conjecture}

\theoremstyle{definition}

\newtheorem{definition}[theorem]{Definition}

\newtheorem{remark}[theorem]{Remark}

\newtheorem{problem}[theorem]{Problem}

\newcommand\E{\mathbb{E}}
\newcommand\Z{\mathbb{Z}}
\newcommand\R{\mathbb{R}}
\newcommand\Q{\mathbb{Q}}

\newcommand\C{\mathbb{C}}

\newcommand\N{\mathbb{N}}

\newcommand\K{\mathcal{K}}

\newcommand\op{{\operatorname{op}}}

\newcommand{\mK}{\mathcal{K}}

\newcommand{\XX}{\mathrm{X}}

\newcommand{\mH}{\mathcal{H}}

\begin{document}


\baselineskip=17pt


\title[]{The structure of non-commutative multiple correlation sequences and applications}

\author[O. Shalom]{Or Shalom}
\address{Department of Mathematics\\ Bar Ilan University \\ 
Ramat Gan \\
5290002, Israel}
\email{Or.Shalom@math.biu.ac.il}

\date{\today}

\begin{abstract}
Let $\Gamma$ be a countable abelian group with Pontryagin dual $\Sigma:=\widehat{\Gamma}$.
A non-commutative $2$-fold multiple correlation sequence is a double-indexed sequence
$$b(\gamma,\gamma') = \int_X T_\gamma S_{\gamma'} f\cdot T_\gamma g\cdot h\,d\mu ,$$
where $f,h\in L^2(\mu)$ and $g\in L^\infty(\mu)$ are functions on a probability space
$\XX=(X,\mathcal{B},\mu)$ and $T,S$ are two, not necessarily commuting, measure preserving
$\Gamma$-actions. We prove that $b$ is of this form if and only if there are two finite
Borel measures $\sigma_1,\sigma_2$ on $\Sigma$ and a bounded operator
$G:L^2(\Sigma,\sigma_1)\rightarrow L^2(\Sigma,\sigma_2)$ with
$$b(\gamma,\gamma')=\int_\Sigma G(\xi_\gamma)\cdot\xi_{\gamma'}\,d\sigma_2 ,$$
where $\xi_\gamma(\chi)=\chi(\gamma)$ is the evaluation map.

We then study the extent to which the triple $(G,\sigma_1,\sigma_2)$ is positive in the
case $f=g=h=1_A$. Positivity fails in the sense available for the spectral measure of a
single correlation sequence, but survives on positive functions with non-negative
Fourier coefficients, where we also obtain a quantitative form by generalizing an inequality of Chu from \cite{Chu}. As an application we prove a
multiple recurrence theorem for products of linear forms, and deduce a simultaneous
partition regularity result for a family of quadratic equations, extending, under a non-degeneracy hypothesis, a theorem of
Frantzikinakis and Host.
\end{abstract}



\maketitle

\section{Introduction}
In \cite{furstenberg1977ergodic} Furstenberg gave an ergodic-theoretic proof for Szemer\'edi's theorem \cite{szemeredi1975sets} on the existence of arbitrarily long arithmetic progressions in sets of positive upper Banach density. The proof introduces the notion of a multiple correlation sequence, which in its simplest form is a sequence \begin{equation}\label{MCS}a(n) = \int_X f_0(x)\cdot \prod_{i=1}^{\ell} f_i(T^{in}x)\,d\mu(x)
\end{equation} where $f_0,\dots,f_{\ell}\in L^\infty(\mu)$ are bounded functions on a probability space $\XX=(X,\mathcal{B},\mu)$ and $T:X\rightarrow X$ is a measure-preserving transformation. More specifically, Furstenberg deduced Szemer\'edi's theorem from the following result.
\begin{theorem}[Multiple recurrence]
    Let $\XX=(X,\mathcal{B},\mu,T)$ be a measure-preserving system and $A\in \mathcal{B}$ is a set of positive measure. Then setting $a(n)$ as in \eqref{MCS} with $f_0=\dots=f_{\ell}=1_{A}$ we have
    $$\liminf_{N\rightarrow\infty}\frac{1}{N} \sum_{n=1}^N a(n)>0.$$
\end{theorem}
Since this work, the Ces\'aro averages of such, and more general, multiple correlation sequences were studied intensively in the literature (see e.g., \cite{cl1,cl2,cl3,fw,bergelson1996polynomial,leibman1998multiple,tao2008norm,walsh2012norm}). The existence of the Ces\'aro limit of \eqref{MCS} was established by Host and Kra \cite{host2005nonconventional} and independently by Ziegler \cite{ziegler2007universal}, and exact formulas are also known \cite{BHK,zieglerformula}. The latter have led to several extensions of Szemer\'edi's theorem (see e.g., \cite{Franpoly,BHK,abs}). Multiple correlation sequences can be generalized beyond \eqref{MCS}. For the sake of this paper we will adapt the following rather general definition.
\begin{definition}[Multiple correlation sequences in general]\label{MCS:def}
Let $\Gamma=(\Gamma,+)$ be a countable abelian group and $\ell\geq 1$. An $\ell$-fold multiple correlation sequence is a sequences of the form
\begin{align*}
    a(\gamma_1,\dots,\gamma_{\ell}) = \int_X f_0(x)\cdot \prod_{i=1}^\ell f_i(T^{(i)}_{ \gamma_i}x)\,d\mu(x)
\end{align*}
where $T^{(1)},\dots,T^{(\ell)}:\Gamma\rightarrow \mathrm{Aut}(X,\mathcal{B},\mu)$ are arbitrary measure-preserving actions of $\Gamma$ on a probability space $(X,\mathcal{B},\mu)$, and $f_0,\dots,f_{\ell}\in L^\infty(\mu)$ are bounded functions.
\end{definition}
Notice how here we define $a$ as a multivariable sequence (instead of setting e.g., $\gamma_1=\dots=\gamma_{\ell}$). This will be crucial in our paper, as opposed to existing literature. This role will be explained below.

Mostly in the case of $\Z$-actions and  commuting transformations $T_1,\dots,T_{\ell}$, these sequences were studied intensively in the literature (see e.g., \cite{MC1,MC3,MCFranHost,MC4,MC5,MC2,Leibman1,Leibman2}). Most of these cited results concern a  decomposition result known as the \emph{nil plus null decomposition}. The most famous result is due to Bergelson, Host and Kra \cite{BHK} regarding the classical multiple correlation sequences.
\begin{theorem}\label{firstnilplusnull}
    Let $\ell>0$ be arbitrary. Then any sequence of the form \eqref{MCS} is the sum of an $\ell$-step nilsequence and a null sequence.
\end{theorem}
The notions nilsequence and null sequence are defined formally in \cite{BHK} and are not needed for this paper. However, we shall now define the former in order to introduce a conjecture of Frantzikinakis (see \cite[Problem 1]{Fran}).
\begin{definition}[Nilsystems, Nilsequences and generalized nilsequences]
Let $k\geq 1$. 
\begin{itemize}
    \item A \emph{$k$-step nilmanifold} is a homogeneous space $X=G/\Lambda$, where $G$ is a $k$-step nilpotent Lie group and $\Lambda\leq G$ is a discrete co-compact subgroup. $X$ is equipped with the quotient topology, the Borel $\sigma$-algebra, and a normalized $G$-invariant measure $\mu=\mu_{G/\Lambda}$. Lastly, every $a\in G$ acts on $X$ by the ($\mu$-preserving) translation $T_a(g\Lambda):=(ag)\Lambda$. Equipped with such an action, $(X,\mathcal{B},\mu,T_a)$ is called a \emph{$k$-step nilsystem}.
    \item A \emph{basic $k$-step nilsequence} is a sequence of the form $n\mapsto F(T_a^n x_0)$, where $(X,\mathcal{B},\mu,T_a)$ is a $k$-step nilsystem, $x_0\in X$, and $F\in C(X)$ is a continuous function.
A \emph{$k$-step nilsequence} is a uniform limit of basic $k$-step nilsequences.
\item Assume furthermore that $X$ is ergodic\footnote{This is not explicitly assumed in Frantzikinakis' original definition, but it is certainly needed, as otherwise any bounded sequence is a generalized nilsequence.}. If one relaxes the condition of continuity and only requires that $F$ is Riemann integrable, that is, bounded and continuous outside of a $\mu$-measure zero set, then a sequence of the form $n\mapsto F(T_a^n x_0)$ is called a \emph{basic generalized $k$-step nilsequence}. Finally, a \emph{generalized $k$-step nilsequence} is a uniform limit of basic generalized $k$-step nilsequences.
\end{itemize}
\end{definition}
All of these definitions can be generalized to arbitrary abelian groups, see e.g., \cite{jt21-1}. Unfortunately, while nilsequences are well understood (see e.g., \cite{gt1}), this is not the case for null-sequences. The latter leads to serious obstacles in applying Theorem \ref{firstnilplusnull} in the setting where one does not have the luxury of taking Ces\'aro averages. In the case of single correlation sequences however, in many applications, see e.g., \cite{FranHost}, it is possible to rely on the Herglotz--Bochner theorem in place of Theorem~\ref{firstnilplusnull}. That is, Herglotz and Bochner provide a full structural classification of single correlation sequences, showing that $$\int_X T_\gamma f\cdot g\,d\mu = \int_{\Sigma} \xi_\gamma~d\sigma_{f,g}$$ for some bounded variation, complex valued measure $\sigma_{f,g}$ on the Pontryagin dual $\Sigma=\widehat{\Gamma}$.\footnote{This was originally established for $\Gamma=\Z$ and $f=g$, but these generalizations are standard, see \cite{folland-harmonic-analysis}.} It is natural to ask whether this result can be generalized to multiple correlation sequences. This is a huge open problem posed by Frantzikinakis in \cite[Problem 1]{Fran}. Here $\Gamma=\mathbb{Z}$, and the actions $T^{(1)},\dots, T^{(\ell)}$ are commutative.
\begin{conjecture}[Frantzikinakis Conjecture]\label{prob:F1}
Let $\ell\geq 1$, let $\XX=(X,\mathcal{B},\mu,T)$ be an ergodic system and let $a(n)$ be an $\ell$-fold \textbf{commutative} multiple correlation sequence. Then, for every $\varepsilon>0$ there exist a complex Borel measure $\nu$ of bounded variation on a compact metric space $S$ and measurable maps $s\mapsto \psi_s(n)\in L^\infty(\nu)$ for all $n\in \mathbb{N}$ from $s\in S$ to a bounded generalized nilsequence $\psi_s(n)$ such that $$\left|a(n) - \int_S \psi_s(n)\,d\nu(s)\right| < \varepsilon$$ for all $n\in \mathbb{N}$. In other words, $a(n)$ can be approximated in $\ell^\infty(\N)$ by integral combinations of generalized nilsequences.
\end{conjecture}
\begin{remark}
In \cite{BrietGreen} Briet and Green proved that the conjecture is false for nilsequences, even if one allows $F$ to be piecewise continuous. In other words the role of \emph{generalized nilsequences} specifically is necessary.
\end{remark}

Our main result concerns a multivariable non-commutative version of this conjecture. However, before we state our main result, we cautiously note that in the non-commutative setting, a result of Frantzikinakis, Lesigne and Wierdl \cite{random} implies that \footnote{This was established in the context of $\Gamma=\Z$ actions, but the argument generalizes naturally to all countable abelian groups.} any bounded sequence is the diagonal of a $2$-fold non-commutative multiple correlation sequence (i.e., taking $\gamma_1=\gamma_2$ in Definition \ref{MCS:def}). Thus, our structural theorem (see Theorem~\ref{main:thm}) is only interesting when $\gamma$ and $\gamma'$ are allowed to differ. In particular, we do not claim to give a satisfactory answer to the commutative problem (see however Problem~\ref{commutative:problem}).\\
 
 Our main result is the following characterization. We call a triple
$(G,\sigma_1,\sigma_2)$ a \emph{spectral triple} over $\Gamma$ if $\sigma_1,\sigma_2$ are
finite Borel measures on $\Sigma:=\widehat{\Gamma}$ and
$G:L^2(\Sigma,\sigma_1)\rightarrow L^2(\Sigma,\sigma_2)$ is a bounded operator. We set
$$b_{(G,\sigma_1,\sigma_2)}(\gamma,\gamma'):=\int_\Sigma G(\xi_\gamma)\cdot\xi_{\gamma'}\,d\sigma_2 ,$$
where $\xi_\gamma:\Sigma\rightarrow S^1$ is the evaluation map $\xi_\gamma(\chi)=\chi(\gamma)$.

\begin{theorem}[Structure theorem for $2$-fold correlation sequences for non-commutative transformations]\label{main:thm}
Let $\Gamma$ be a countable abelian group and let $b:\Gamma\times\Gamma\rightarrow\C$. The
following are equivalent.
\begin{itemize}
\item[(i)] $b$ is a $2$-fold (non-commutative) multiple correlation sequence. Namely, there
are a probability space $\XX=(X,\mathcal{B},\mu)$, two (not necessarily commuting)
$\Gamma$-actions $T,S:\Gamma\rightarrow\mathrm{Aut}(\XX)$ and functions $f,h\in L^2(\mu)$,
$g\in L^\infty(\mu)$ with
\begin{equation}\label{strthm}
b(\gamma,\gamma')=\int_{\XX} T_\gamma S_{\gamma'} f \cdot T_\gamma g \cdot h ~d\mu.
\end{equation}
\item[(ii)] $b=b_{(G,\sigma_1,\sigma_2)}$ for some spectral triple $(G,\sigma_1,\sigma_2)$.
\end{itemize}
Moreover, in \emph{(i)}$\Rightarrow$\emph{(ii)} one may take $\sigma_1$ to be the spectral
measure of $\overline{h}$ with respect to $T$ and $\sigma_2$ that of $f$ with respect to
$S$, and then $\|G\|_{\op}\leq\|g\|_{L^\infty(\mu)}$, and in
\emph{(ii)}$\Rightarrow$\emph{(i)} one may take $\XX$ to be a Gaussian system on a standard
probability space and $g\equiv\mathbf 1$.
\end{theorem}

The first implication is proved in Section \ref{firstdir} and the second in Section
\ref{seconddir}. In Section \ref{positivity} we study to what extent the triple $(G,\sigma_1,\sigma_2)$ is
positive when $f=g=h=1_A$, and in Section \ref{prapp} we use the above
positivity results to prove a multiple recurrence theorem for products of
linear forms, and deduce from it a simultaneous partition regularity result.

It is natural to ask whether this result could be applied to resolve Frantzikinakis' conjecture. We believe that this might be possible if one is able to solve the following classification problem.
\begin{problem}\label{commutative:problem}
Determine for which triples $(G,\sigma_1,\sigma_2)$ the expression on the right hand side of \eqref{strthm} is a multiple correlation sequence for commuting $T$ and $S$.
\end{problem}
The same problem in the case where $T=S$ is also interesting.
\subsection*{Acknowledgments}
This research was supported by Alon Fellowship, an NSF grant DMS-1926686 and an ISF grant 3056/21. Thanks to Claude Opus I was able to find a mistake in an earlier version of this manuscript \cite{ShalomGrothendieck} (see Proposition~\ref{posfail}). The fix necessitated proving a quantitative lower bound in the spirit of \cite[Lemma 1.6]{Chu} for multicorrelation sequences associated with PSD doubly-stochastic operators (see Theorem~\ref{conj}). Claude Opus was also used for copy-editing; see in particular Remark~\ref{inconsistent}.
\section{The first direction: spectral measures}\label{firstdir}
In this section we prove the first direction of Theorem~\ref{main:thm}. Throughout, we allow $\Gamma$
to be an arbitrary countable abelian group, but the result is already new for
$\Gamma=\mathbb{Z}$.

Let $\Gamma$ be a countable abelian group and let $T,S:\Gamma\rightarrow\mathrm{Aut}(\XX)$
be two (not necessarily commuting) $\Gamma$-actions on a probability space
$\XX=(X,\mathcal{B},\mu)$. By standard abuse of notation we denote by  $T,S:\Gamma\rightarrow\mathcal{U}(L^2(\mu))$ the corresponding unitary
representations. Let $\Sigma$ denote the Pontryagin dual of $\Gamma$. Gelfand theory gives
rise to a $\star$-homomorphism $C(\Sigma)\rightarrow\mathcal{L}(L^2(\XX))$ sending a
continuous $\phi:\Sigma\rightarrow\mathbb{C}$ to an operator $T_\phi$ with
$\|T_\phi\|_{\op}\leq\|\phi\|_\infty$, determined by $T_{\xi_\gamma}=T_\gamma$.
We write $S_\psi$ for the operators attached in the same way to the representation $S$.

\subsection{Spectral measures}
For a unitary representation $T$ of $\Gamma$ on $L^2(\mu)$ and $u\in L^2(\mu)$,
the Bochner--Herglotz theorem (see e.g., \cite{folland-harmonic-analysis}) implies the existence of a unique positive finite Borel measure
$\sigma^{T}_{u}$ on $\Sigma$, called the \emph{spectral measure of $u$ with respect to $T$}, such that
$$\langle T_\gamma u,u\rangle=\int_\Sigma\xi_\gamma\,d\sigma^{T}_u,$$ for all $\gamma\in \Gamma$.
The following lemma is classical, see e.g., \cite{folland-harmonic-analysis}.
\begin{lemma}\label{specid}
For every $u\in L^2(\mu)$ and every $\eta\in C(\Sigma)$ we have
$\langle T_\eta u,u\rangle=\int_\Sigma\eta\,d\sigma^{T}_u$. Consequently,
$$\sigma^{T}_u(\Sigma)=\|u\|_{2}^{2}\qquad\text{and}\qquad
\big\|T_\phi u\big\|_{L^2(\mu)}^{2}= \|\phi\|_{L^2(\sigma^{T}_u)}^2,$$ for all $\phi\in C(\Sigma)$.
\end{lemma}

\subsection{The bilinear form}
Fix $f,h\in L^2(\XX)$ and $g\in L^\infty(\XX)$. For $\phi,\psi\in C(\Sigma)$ the function
$T_\phi(S_\psi f\cdot g)$ lies in $L^2(\XX)$, so
$$\Phi(\phi,\psi):=\int_X T_\phi\big(S_\psi f\cdot g\big)\cdot h\,d\mu$$
is well defined and bilinear. The following estimate is the key point.

\begin{proposition}\label{keybound}
Set $\sigma_1:=\sigma^{T}_{\overline h}$ and $\sigma_2:=\sigma^{S}_{f}$. Then for all $\phi,\psi\in C(\Sigma)$, we have
$$\big|\Phi(\phi,\psi)\big|\ \leq\ \|g\|_{L^\infty(\mu)}\,
\|\phi\|_{L^{2}(\sigma_1)}\,\|\psi\|_{L^{2}(\sigma_2)}.$$
\end{proposition}

\begin{proof}
Write $w:=S_\psi f\cdot g$. Since $\int_X u\,h\,d\mu=\langle u,\overline h\rangle$ and
$T_\phi^{*}=T_{\overline\phi}$,
$$\Phi(\phi,\psi)=\big\langle T_\phi w,\ \overline h\big\rangle
=\big\langle w,\ T_{\overline\phi}\,\overline h\big\rangle .$$
By the Cauchy--Schwarz inequality and Lemma \ref{specid},
$$|\Phi(\phi,\psi)|\leq\|w\|_2\,\big\|T_{\overline\phi}\overline h\big\|_2
\leq\|g\|_\infty\|S_\psi f\|_2\,\big\|T_{\overline\phi}\overline h\big\|_2
=\|g\|_\infty\,\|\psi\|_{L^2(\sigma^{S}_f)}\,\|\phi\|_{L^2(\sigma^{T}_{\overline h})}.$$
\end{proof}

\begin{theorem}[First direction of Theorem \ref{main:thm}]\label{firstdirthm}
With $\sigma_1,\sigma_2$ as in Proposition \ref{keybound} there is a bounded operator
$$G:L^{2}(\Sigma,\sigma_1)\longrightarrow L^{2}(\Sigma,\sigma_2),
\qquad \|G\|_{\op}\leq\|g\|_{L^\infty(\mu)},$$
such that
\begin{equation}\label{spec2gen}
\Phi(\phi,\psi)=\int_X T_\phi\big(S_\psi f\cdot g\big)\cdot h\,d\mu
=\int_\Sigma G(\phi)\cdot\psi\,d\sigma_2,
\end{equation}
for all $\phi,\psi\in C(\Sigma).$
In particular, specialising to $\phi=\xi_\gamma$ and $\psi=\xi_{\gamma'}$,
\begin{equation}\label{spec2}
\int_X T_{\gamma}S_{\gamma'}f\cdot T_\gamma g\cdot h~d\mu
=\int_\Sigma G(\xi_\gamma)\cdot\xi_{\gamma'}~d\sigma_2 .
\end{equation}
\end{theorem}

\begin{proof}
By Proposition \ref{keybound} the form $\Phi$ is bounded for the norms of
$L^2(\sigma_1)\times L^2(\sigma_2)$. In particular it vanishes whenever one of its
arguments is null for the corresponding measure, so it descends to the images of
$C(\Sigma)$ in these spaces. Since this image is dense, we conclude that $\Phi$ extends uniquely to a bounded bilinear
form on $L^2(\sigma_1)\times L^2(\sigma_2)$ with the same bound. The map
$(\phi,\psi)\mapsto\Phi(\phi,\overline\psi)$ is then a bounded sesquilinear form, so by the
Riesz representation theorem there is a bounded
$G:L^2(\sigma_1)\rightarrow L^2(\sigma_2)$ with $\|G\|_{\op}\leq\|g\|_\infty$ and
$\Phi(\phi,\overline\psi)=\langle G\phi,\psi\rangle_{L^2(\sigma_2)}$, that is
$\Phi(\phi,\psi)=\int_\Sigma G(\phi)\psi\,d\sigma_2$.
\end{proof}

\section{The second direction of Theorem~\ref{main:thm}}\label{seconddir}
Our goal in this section is to show that any spectral triple corresponds to a genuine $2$-fold multiple correlation sequence. In this section we will work with vector spaces defined both over the reals $\R$ and over the complex numbers $\C$. To distinguish the two, we will often add a subscript, $V_{\R}$ and $V_{\C}$ for vector spaces over the reals and complex numbers, respectively.

Throughout we fix $\Sigma:=\widehat{\Gamma}$, two finite Borel measures $\sigma_1,\sigma_2$
on $\Sigma$, and a bounded operator
$G:L^2(\Sigma,\sigma_1)\rightarrow L^2(\Sigma,\sigma_2)$ and set
$$b_{(G,\sigma_1,\sigma_2)}(\gamma,\gamma') = \int_{\Sigma} G(\xi_\gamma)\cdot \xi_{\gamma'}~d\sigma_2 .$$
As before, $\xi_\gamma:\Sigma\rightarrow S^1$ denotes the evaluation map by $\gamma$. For
$i=1,2$ we denote by $M^{(i)}_{\xi_\gamma}$ the operator of multiplication by $\xi_\gamma$
on $L^2_{\C}(\sigma_i)$, which is clearly unitary.

Since $G$ maps between two different spaces, we first reduce to a single one. Set
$$\mH:=L^2_{\C}(\Sigma,\sigma_1)\oplus L^2_{\C}(\Sigma,\sigma_2),\qquad
M_{\xi_\gamma}:=M^{(1)}_{\xi_\gamma}\oplus M^{(2)}_{\xi_\gamma},\qquad
\widetilde G(u_1,u_2):=(0,Gu_1),$$
and
$$x_0:=(\bm{1}_{\Sigma},0),\qquad y_0:=(0,\bm{1}_{\Sigma}),$$ where $\bm{1}_{\Sigma}:\Sigma\rightarrow \C$ is identically $1$.
Then $\mH$ is a separable complex Hilbert space, $\gamma\mapsto M_{\xi_\gamma}$ is a unitary
representation of $\Gamma$ on $\mH$, $\widetilde G$ is bounded with
$\|\widetilde G\|_{\op}=\|G\|_{\op}$, and $\|x_0\|^2=\sigma_1(\Sigma)$,
$\|y_0\|^2=\sigma_2(\Sigma)$. The following observation is key:
\begin{equation}\label{bGl}
    b_{(G,\sigma_1,\sigma_2)}(\gamma,\gamma')
    = \langle M_{\xi_{\gamma'}}\,\widetilde G\, M_{\xi_{\gamma}}x_0,\ y_0\rangle_{\mH}.
\end{equation}
Indeed,
$M_{\xi_\gamma}x_0=(\xi_\gamma,0)$, so $\widetilde GM_{\xi_\gamma}x_0=(0,G\xi_\gamma)$ and
$M_{\xi_{\gamma'}}\widetilde GM_{\xi_\gamma}x_0=(0,\xi_{\gamma'}G(\xi_\gamma))$, whose inner
product with $y_0=(0,\mathbf 1_{\Sigma})$ is $\int_\Sigma G(\xi_\gamma)\xi_{\gamma'}\,d\sigma_2$.

From this point on the argument depends only on the data
$(\mH,M,\widetilde G,x_0,y_0)$, and not on the particular form of $\mH$. If $G=0$, then $b_{(G,\sigma_1,\sigma_2)}=0$ and is therefore a multiple correlation sequence with (say) $f=0$. Thus, we shall assume throughout that $\|G\|_{op}\not = 0$. Furthermore, we may normalize (say by multiplying $f$ by a constant later) and assume that $\|G\|_{op}=1$.

We see from \eqref{bGl} that $b_{(G,\sigma_1,\sigma_2)}$ is an inner product involving two $\Gamma$-representations, separated by the operator $G$. For technical reasons, this separation is problematic for us. Our first goal is to find an alternative representation for $b_{(G,\sigma_1,\sigma_2)}$ that involves the inner product of two $\Gamma$-representations without a separating operator. This will require replacing $L^2_{\C}(\sigma_1)$ and $L^2_{\C}(\sigma_2)$ with a different space. Unfortunately, we know nothing about $G$ other than the fact that $\|G\|_{\op} = 1$. As a first step we shall use that and the Halmos-dilation \cite{Halmosdilation} in order to replace $G$ with a unitary operator. The formal argument is given below.

\subsection{The Halmos dilation}
Throughout this section, we let $\mH=\mH_{\C}$ denote a complex, separable Hilbert space. A \emph{contraction} on $\mathcal{H}$ is a bounded linear operator $C:\mathcal{H}\rightarrow \mathcal{H}$ with $\|C\|_{op}\leq 1$. Given such a contraction, the operators $I-CC^*$ and $I-C^*C$ are self-adjoint and non-negative and as such they admit square roots. In fact, we can choose this root canonically by fixing $0<c_n\in\R$ such that $\sqrt{1-t} = 1-\sum_{n\geq 1} c_nt^n$ (for all $|t|\leq 1$), and noting that for a self-adjoint operator $P$ with $\|P\|_{op}\leq 1$, the operator $R(P) = I-\sum_{n\geq 1} c_n P^n$ is a well defined self-adjoint operator such that $R(P)^2=I-P$. With the notations described above, we can now define the Halmos dilation.
\begin{definition}[Halmos dilation]\label{Halmosdilation:def}
    Let $C:\mathcal{H}\rightarrow \mathcal{H}$ be a contraction, and set $D:=R(C^*C)$ and $D_*:=R(CC^*).$ The Halmos dilation of $C$ is the operator $\widehat{C}$ defined on $\mathcal{K}:=\mathcal{H}\oplus \mathcal{H}$ by
    $$\widehat{C} = \begin{pmatrix}C & D_* \\ D & -C^*\end{pmatrix}.$$
\end{definition}
The following properties of the Halmos dilation are classical. Crucially, property $(iii)$ is the desired unitary.
\begin{lemma}\label{HDproperties}
    In the setting of Definition~\ref{Halmosdilation:def} we have
    \begin{itemize}
        \item[(i)] $D^2=I-C^*C$ and $D_*^2=I-CC^*$.
        \item[(ii)] $CD=D_*C$ and $C^*D_*=DC^*$. 
        \item[(iii)] $\widehat{C}$ is unitary.
        \item[(iv)] for all $a,b\in \mathcal{H}$, $\langle\widehat{C}(a,0),(b,0)\rangle_{\mathcal{K}} = \langle Ca,b\rangle_{\mathcal{H}}.$
    \end{itemize}
\end{lemma}
\begin{proof}
    Property $(i)$ is immediate from the definition of $R$. Moreover, since $R$ is a power series and $C(C^*C) = (CC^*)C$, opening the brackets gives $(ii)$. To prove $(iii)$ observe that
    $$\left(\widehat{C}\right)^* = \begin{pmatrix}C^* & D \\ D_* & -C\end{pmatrix}.$$ Multiplying with $\widehat{C}$ we get
    $$\left(\widehat{C}\right)^*\widehat{C} =\begin{pmatrix}C^*C+D^2 & C^*D_*-DC^* \\ D_*C-CD & D_*^2+CC^*\end{pmatrix}.$$
    Since $D^2 = I-C^*C$ and $D_*^2 = I-CC^*$, the diagonal trivializes. Furthermore, from $(ii)$ we see that the off-diagonal is zero. The same type of computation also shows that $\widehat{C}\left(\widehat{C}\right)^*=I.$ Finally $(iv)$ is a direct computation.
\end{proof}

Now, let $\mathcal{K}=\mathcal{H}\oplus \mathcal{H}$ (here $\mK$ is a complex vector space). Let $M^{\mK}_{\xi_\gamma} = M_{\xi_{\gamma}}\times M_{\xi_{\gamma}}$, and let $\widehat{G}$ denote the Halmos dilation for the contraction $\widetilde G$ on $\mH$. Finally, let $U_{\gamma}:=\widehat{G}\circ M^{\mK}_{\xi_\gamma}\circ \widehat{G}^{-1}.$ Then taking $x=\widehat{G}(x_0,0)$ and $y=(y_0,0)$, we have
\begin{equation}\label{twounitaries}
    b_{(G,\sigma_1,\sigma_2)}(\gamma,\gamma') = \langle M^{\mK}_{\xi_{\gamma'}} U_\gamma x,y\rangle_{\mK}. 
\end{equation}
Indeed, $$U_\gamma x = \widehat{G}\circ M_\gamma^{\mK} \circ \widehat{G}^{-1} \circ \widehat{G}(x_0,0)= \widehat{G} (M_{\xi_\gamma}x_0,0),$$ and so, using that $M^{\mK}_{\xi_{\gamma'}}$ is unitary with adjoint $M^{\mK}_{\xi_{-\gamma'}}$,
\begin{align*}\langle M^{\mK}_{\xi_{\gamma'}} U_\gamma x,y\rangle_{\mK} &= \langle M^{\mK}_{\xi_{\gamma'}} \widehat{G}(M_{\xi_\gamma}x_0,0),(y_0,0)\rangle_{\mK}\\&= \langle \widehat{G} (M_{\xi_\gamma}x_0,0),(M_{\xi_{-\gamma'}}y_0,0)\rangle_{\mK} \\&= \langle \widetilde G\,M_{\xi_\gamma}x_0,\ M_{\xi_{-\gamma'}}y_0\rangle_{\mH}\\&= \langle M_{\xi_{\gamma'}}\widetilde G\,M_{\xi_\gamma}x_0,\ y_0\rangle_{\mH}\\&=b_{(G,\sigma_1,\sigma_2)}(\gamma,\gamma'),
\end{align*}
where the third equality follows from Lemma~\ref{HDproperties}(iv) and the last equality from \eqref{bGl}.

Now that there are only two unitary representations involved in \eqref{twounitaries}, with no separating operator, we are set to construct our probability space.

\subsection{The Gaussian system}
The construction below is the classical Gaussian functor, which associates to a real Hilbert
space and an orthogonal representation on it a measure preserving system; see
\cite{glasner2015ergodic} and \cite{Svante} for detailed treatments. We include the details
we need since we shall use the explicit form of the correspondence. For technical reasons, since the Gaussian measure is defined on $\R$, rather than $\C$, we must initially work with a real Hilbert space $\mH_{\R}$ (later we take $\mH_{\R}$ denote the real Hilbert space underlying $\mK$ from the previous section, see Section~\ref{complex}). We also fix throughout a countable subgroup $\Lambda$ of the orthogonal group of $\mH_{\R}$, and choose a dense and countable subset $D_0\subseteq \mH_{\R}$ (the choice of $D_0$ will not matter). Then the set
$$
D:=\Big\{\textstyle\sum_{j=1}^{n}q_j A_j(v_j)\ :\ n\ge0, q_j\in\Q, A_j\in\Lambda, v_j\in D_0 \Big\}
$$
is a $\Lambda$-invariant, countable, $\Q$-linear, dense subspace of $\mH_{\R}$. With respect to this $D$, we set $X=\R^D$ and equip $X$ with the Borel $\sigma$-algebra $\mathcal{B}$. Since $D$ is countable, $X$ is a Polish space and the $\sigma$-algebra is generated by the coordinate maps
\begin{align*}
    \widehat{v}:X\rightarrow \R\\
    \widehat{v}(w) =w_v
\end{align*} where $v\in D$. It follows from Kolmogorov extension theorem that there exists a unique measure $\mu$ on $X$ with the property that for every finite tuple $v_1,\dots,v_n$ of distinct elements in $D$, the push-forward of $\mu$ under $\widehat{v}_1\times\cdots\times \widehat{v}_n$ is the Gaussian measure on $\R^n$ with covariance matrix $\left(\langle v_i,v_j\rangle\right)_{i,j=1}^n$. We refer to $\XX = (X,\mathcal{B},\mu)$ as the \emph{Gaussian probability space over $\mathcal{H}_{\R}$} (relative to $D$).
The group $\Lambda$ acts on $X$. Indeed, for any $A\in \Lambda$ we may define an action
\begin{equation}\label{R}(R_A x)_v = x_{A^{-1}v}.
\end{equation}
This is well defined because $D$ is $\Lambda$-invariant, and one checks directly from \eqref{R} that $R_A\circ R_B = R_{AB}$ and $R_{\mathrm{id}}=\mathrm{id}$, so each $R_A$ is a bijection with inverse $R_{A^{-1}}$. Furthermore, $R_A$ is continuous (and therefore measurable) and each $R_A$ preserves $\mu$. Indeed, let $v_1,\dots,v_n\in D$ be distinct. By \eqref{R}, we have $\widehat{v_i}\circ R_A = \widehat{A^{-1}v_i}$, and $A^{-1}v_1,\dots,A^{-1}v_n$ are distinct elements of $D$. Hence, the push-forward of $\mu$ under $(\widehat{v}_1\times\cdots\times\widehat{v}_n)\circ R_A$ is the Gaussian measure on $\R^n$ with covariance matrix $$\left(\langle A^{-1}v_i,A^{-1}v_j\rangle\right)_{i,j=1}^n = \left(\langle v_i,v_j\rangle\right)_{i,j=1}^n,$$ the equality holding because $A$ is orthogonal. This is precisely the push-forward of $\mu$ under $\widehat{v}_1\times\cdots\times\widehat{v}_n$, so $(R_A)_*\mu$ and $\mu$ agree on all cylinder sets. As these generate $\mathcal{B}$, we conclude that $(R_A)_*\mu=\mu$. We refer to $(X,\mathcal{B},\mu,\Lambda)$ as the \emph{Gaussian $\Lambda$-system over $\mH$} (relative to $D$).
\subsubsection{Complexification}\label{complex}
Let $\mK_{\R}$ denote the underlying real vector space of $\mK$. Namely, as sets (and as abelian groups under addition), $\mK_\R=\mK$, and the inner product is given by
$$\langle u,v\rangle_{\mK_{\R}} = \mathrm{Re}\langle u,v\rangle_{\mK}.$$ This turns $\mK_\R$ into a real separable Hilbert space. Furthermore, if $(f_j)_{j\in \N}$ is an orthonormal basis for $\K$ over $\C$, then the family $(f_j)_{j\in \N}$ and $(i\cdot f_j)_{j\in \mathbb{N}}$ where $i=\sqrt{-1}$, is an orthonormal basis for $\K_{\R}$ over $\R$. Now, let $\K_{\C} = \K_{\R}\otimes_{\R} \C$ denote the algebraic tensor product (i.e., the  completion of the set of all finite formal sums of tensors $f\otimes t$ where $f\in \K_{\R}$ and $t\in \C$, where addition is bilinear). This is now a vector space over $\C$, where the scalar-multiplication by $s\in \C$ and inner product are given by
$$s\cdot (f\otimes t) = f\otimes (st), \qquad \langle f\otimes t, g\otimes s\rangle_{\K_{\C}} = t\overline{s}\langle f,g\rangle_{\K_{\R}}.$$
Note that since $\C = \R\oplus i\R$, the sums $$f\otimes 1+g\otimes i$$ for $f,g\in \K_\R$ generate $\K_{\R}\otimes_{\R}\C$. 

Generally, $\K_{\C}$ is bigger than $\K$ (technically $\dim \mK_{\C} = 2\dim \mK$, however generally both are infinite dimensional). Yet, we may embed $\K$ in $\K_{\C}$ by setting:
\begin{align*}\imath:&\mK \rightarrow \K_{\C}\\
\imath(f) &= \frac{1}{\sqrt{2}}\left(f\otimes 1 -Jf\otimes i\right),
\end{align*}
where $Jf=i\cdot f$. Unlike the natural embedding $f\mapsto f\otimes 1$, this embedding is $\C$-linear. Indeed, since $J^2f = -f$ we have
\begin{align*}
    \imath(Jf) &= \frac{1}{\sqrt{2}} (Jf\otimes 1 - \left(J^2 f)\otimes i\right)\\&= \frac{1}{\sqrt{2}} (Jf\otimes 1 +f\otimes i)\\&= \frac{1}{\sqrt{2}}(i\cdot f\otimes 1 - i\cdot (Jf\otimes i))\\&= i\cdot \imath(f).
\end{align*}

\begin{lemma}\label{commute}
    Let $A\in O(\K_{\R})$ be an orthogonal operator. The map 
    $$A_{\C} (f\otimes s) = Af\otimes s$$ extends linearly to a unitary on $\K_{\C}$. Furthermore, if $A$ is $\C$-linear, then  $A_{\C}\circ \imath = \imath\circ A.$
\end{lemma}
\begin{proof}
    We have 
    $$\langle A_{\C} (f\otimes s),A_{\C} (f\otimes s)\rangle = |s|^2 \langle Af,Af\rangle = |s|^2\langle f,f\rangle  = \langle f\otimes s,f\otimes s\rangle.$$ Hence, $A_{\C}$ preserves the inner product on elementary tensors, and hence on all of $\K_{\C}$ by linearity. It is therefore unitary. The assertion $A_{\C}\circ \imath = \imath\circ A$ follows immediately from the assumption that $A$ commutes with $J$.
\end{proof}
\subsubsection{The measure-preserving system}
Since every unitary on $\K$ is an orthogonal matrix on $\mK_{\R}$ (indeed, it preserves $\langle\cdot,\cdot\rangle_{\mK_{\R}}$ because $\langle Au,Av\rangle_{\mK_{\R}}=\mathrm{Re}\langle Au,Av\rangle_{\mK}=\mathrm{Re}\langle u,v\rangle_{\mK}=\langle u,v\rangle_{\mK_{\R}}$), we may let $\Lambda$ denote the group generated by $M_{\xi_{\gamma}}^{\mK}$ for all $\gamma\in \Gamma$ and $\widehat{G}$, the Halmos dilation of $\widetilde G$ (so it also contains $\{U_{\gamma} : \gamma\in \Gamma\}$). Fix $D\subseteq \mK_{\R}$ as above and let $\XX=(X,\mathcal{B},\mu,\Lambda)$ be a Gaussian $\Lambda$-system over $\mK_{\R}$ relative to $D$. For $A\in O(\mK_{\R})$, we write $R_A$ for the measure-preserving transformation associated with $A$, as defined in \eqref{R}. 
\begin{lemma}
    There exists an isometry $\Psi:\mK\rightarrow L^2(\mu)$ such that $$\Psi(A^{-1} v)=\Psi(v)\circ R_A$$ for all $v\in \mK$ and all $A\in \Lambda$. 
\end{lemma}
\begin{proof}
Recall that for $v\in D$ we have a map $\widehat{v}:X\rightarrow \R$ defined by $\widehat{v}(w) = w_v$. The definition of $\mu$ implies that $v\mapsto \widehat{v}\in L^2(\mu)$ is an isometry on $D$. Indeed,
\begin{align*} \|\widehat{v}-\widehat{u}\|_{L^2(\mu)}^2 &= \mathbb{E}_{\mu} \widehat{u}^2 - 2\mathbb{E}_{\mu} \widehat{u}\widehat{v} + \mathbb{E}_\mu \widehat{v}^2  \\&= \|u\|_{\mK_\R}^2 - 2\langle u,v\rangle_{\mK_\R} + \|v\|_{\mK_\R}^2 = \|u-v\|_{ \mK_\R}^2,  
\end{align*}
where the second equality follows from the assumption on the covariance matrices which determines $\mu$. A similar computation also shows that $v\mapsto \widehat{v}$ is linear as $\|\widehat{sv+tu} - s\widehat{v}- t\widehat{u}\|_{L^2(\mu)}=0$ for all $u,v\in D$ and $s,t\in \Q$. Thus, since $D$ is dense and $L^2(\mu)$ is complete, $v\mapsto\widehat{v}$ extends to an $\R$-linear isometry on $\mK_\R$, and furthermore $\E_{\mu} \widehat{u}\widehat{v} = \langle u,v\rangle$, now for all $u,v\in \K_{\R}$. We may therefore let $W:\K_\R\otimes_\R\C\longrightarrow L^2_\C(\mu)$ denote the linear map determined by $$W\big(a\otimes1+b\otimes i\big):=\widehat a+i\,\widehat b, \qquad a,b\in \mK_{\R}.$$ Note that $W$ is an isometry. Indeed, since $\widehat c,\widehat d$ are real valued,
$$\langle W(a\otimes1+b\otimes i),W(c\otimes1+d\otimes i)\rangle_{L^2(\mu)}
=\mathbb{E}_\mu\big[(\widehat a+i\widehat b)(\widehat c-i\widehat d)\big]
=\langle a,c\rangle+\langle b,d\rangle+i\big(\langle b,c\rangle-\langle a,d\rangle\big),$$
which is precisely $\langle a\otimes1+b\otimes i,\,c\otimes1+d\otimes i\rangle_{\K_\C}$. Now, set $\Psi = W\circ \imath:\mK\rightarrow L^2_{\C}(\mu)$. Since $W$ and $\imath$ are isometries, so is $\Psi$. Since $A^{-1}$ is $\C$-linear, Lemma~\ref{commute} gives
\begin{align*}\Psi(A^{-1} v) = W\circ \imath (A^{-1}v)=W \circ A^{-1}_{\C}\circ \imath(v). 
\end{align*}
It suffices to show that $R_A\circ W = W\circ A^{-1}_{\C}$. To see this observe that for all $a,b\in D$ we have,
$$R_A (\widehat{a} + i\widehat{b}) = R_A (\widehat{a}) + iR_A(\widehat{b}) = \widehat{A^{-1}a} + i\widehat{A^{-1}b}=W\big(A^{-1}a\otimes1+A^{-1}b\otimes i\big)=W\big(A^{-1}_\C(a\otimes1+b\otimes i)\big).$$ Since $R_A$ is measure-preserving, by taking limits this equality extends to all $a,b\in \mK_\R$ and the proof is now complete.
\end{proof}
\subsection{Concluding the proof}
We finally put everything together. From \eqref{twounitaries} we have
$$b_{(G,\sigma_1,\sigma_2)}(\gamma,\gamma') = \langle M_{\xi_{\gamma'}}^{\mK} U_\gamma x,y\rangle$$ for $x=\widehat{G}(x_0,0)$ and $y=(y_0,0)$, both in $\mK$. Fix $f=\Psi(x)$, $g=1$ and $h=\overline{\Psi(y)}$ and set $T_\gamma:= R^{-1}_{M^{\mK}_{\xi_{\gamma}}}$ and $S_\gamma :=R^{-1}_{U_{\gamma}} $. We have,
\begin{align*}
\int_X T_{\gamma'} S_{\gamma}f\cdot T_{\gamma'} g\cdot h\,d\mu
&=\int_X \Psi\big(M_{\xi_{\gamma'}}^{\mK} U_\gamma x\big)\cdot\overline{\Psi(y)}\,d\mu
\\&=\big\langle \Psi\big(M_{\xi_{\gamma'}}^{\mK} U_\gamma x\big),{\Psi(y)}\big\rangle_{L^2(\mu)}
\\&=\big\langle M_{\xi_{\gamma'}}^{\mK}U_\gamma x,y\big\rangle_\K
\\&=b_{(G,\sigma_1,\sigma_2)}(\gamma,\gamma').
\end{align*}
Finally, recall that we normalized $\|G\|_{op}=1$. For a general $G\neq0$ we apply the above to $G/\|G\|_{op}$ and replace $f$ by $\|G\|_{op}\cdot f$, which multiplies the left-hand side by $\|G\|_{op}$ and yields $b_{(G,\sigma_1,\sigma_2)}$.
Note that here $T$ acts via $\gamma'$ and $S$ via $\gamma$, which is the opposite of how we
initially defined it. This is not an issue. Indeed, since $T_{\gamma'}$ and $S_\gamma$ are
measure preserving, we have
$$\int_{\XX} T_{\gamma'} S_\gamma f\cdot h\, d\mu = \int_{X} f\cdot S_{-\gamma} T_{-\gamma'}h\,d\mu.$$
Therefore, setting
$$\widetilde{T}_\gamma := S_{-\gamma},\qquad \widetilde{S}_{\gamma'} := T_{-\gamma'},\qquad
\widetilde{f}:=h,\qquad \widetilde{g}:=\bm{1},\qquad \widetilde{h}:=f,$$
we obtain
$$\int_{\XX} \widetilde{T}_{\gamma}\widetilde{S}_{\gamma'} \widetilde{f}\cdot
\widetilde{T}_{\gamma}\widetilde{g}\cdot \widetilde{h}\,d\mu = b_{(G,\sigma_1,\sigma_2)}(\gamma,\gamma'),$$
which is of the required form. The proof is now complete.
\section{On positivity of the spectral triple \texorpdfstring{$(G,\lambda,\lambda)$}{(G,lambda,lambda)}}\label{positivity}
In this section we isolate a question of independent interest. Suppose that $f=g=h=1_A$ and
$S=T$. To what extent is the triple produced by Theorem \ref{main:thm}
positive? Unfortunately, the situation here is completely different from the case of
spectral measures associated with single correlation sequences. Yet, we are still able to
recover some hidden positivity results, which are later used to obtain a combinatorial
application.

In this situation the two measures of Theorem \ref{firstdirthm} coincide: since
$\overline{1_A}=1_A$ and $S=T$, we have
$$\sigma_1=\sigma^{T}_{\overline{h}}=\sigma^{T}_{1_A}=\sigma^{S}_{f}=\sigma_2 ,$$
a single measure of total mass $\|1_A\|_2^{2}=\delta$, where $\delta:=\mu(A)$. We therefore
set
$$\lambda:=\frac{1}{\delta}\,\sigma^{T}_{1_A},$$
a Borel \emph{probability} measure on $\Sigma$, and rescale accordingly: There is a bounded operator $G$ on $L^2(\Sigma,\lambda)$ with
\begin{equation}\label{Gbound}
\|G\|_{\op}\ \leq\ 1
\end{equation}
such that
\begin{equation}\label{Phiform}
\Phi(\phi,\psi):=\int_\Sigma G(\phi)\cdot\psi\,d\lambda
=\int_X T_\phi\big(T_\psi f\cdot f\big)\cdot f\,d\mu
\end{equation}
for all $\phi,\psi\in C(\Sigma)$.

Now suppose that the Fourier coefficients of $\phi$ are symmetric (i.e., if $\phi = \sum_{\gamma} c_\gamma \xi_\gamma$, then $c_\gamma=c_{-\gamma}$). Then, \eqref{Phiform} may be rewritten as
\begin{equation}\label{Phisym}
\Phi(\phi,\psi)=\int_X f\cdot\big(T_{\phi} f\big)\cdot\big(T_\psi f\big)\,d\mu ,
\end{equation}
which exhibits the symmetry $\Phi(\phi,\psi)=\Phi(\psi,\phi)$ whenever $\psi$ also has symmetric Fourier coefficients.

\subsection{Comparison with the spectral measure}
For two functions rather than three the situation is completely different. It is well known that in this case the spectral measure $\sigma_f$ is positive. In particular, if $\phi\geq 0$, this guarantees that 
$$\int_X T_\phi f\cdot f\,d\mu=\int_\Sigma \phi\,d\sigma_f\geq 0 .$$
In particular, in this scenario one may discard part of the domain, for instance the following inequality was used by Frantzikinakis and Host \cite{FranHost}
$\int_\Sigma\phi\,d\sigma_f\geq\phi(\chi_{\mathrm{triv}})\,\sigma_f(\{\chi_{\mathrm{triv}}\})$ for $\phi\geq0$, where $\chi_{\mathrm{triv}}\in \Sigma$ is the trivial element. Unfortunately, neither of these properties holds for the spectral triplet $(G,\lambda,\lambda)$. The following proposition shows that.

\begin{proposition}[Positivity fails]\label{posfail}
There exist a finite $\Gamma$-system and a set $A$, and functions $\phi,\psi\in C(\Sigma)$
with $\phi\geq0$ and $\psi\geq0$ pointwise on $\Sigma$, such that $\Phi(\phi,\psi)<0$.
Consequently, $G$ need not map non-negative functions to non-negative functions, and the
bilinear form $\Phi$ is not non-negative even when $\phi,\psi$ are. 
\end{proposition}

\begin{proof}
Let $\Gamma=\mathbb{Z}/4\mathbb{Z}$ act on $X=\mathbb{Z}/4\mathbb{Z}$ with the uniform
measure by translation, and let $A=\{0,1,3\}$, so that $\delta=\tfrac34$ and
$f=1_A$. Identify $\Sigma=\widehat{\Gamma}$ with $\mathbb{Z}/4\mathbb{Z}$ via
$\chi_k(t)=i^{kt}$, and let
$$\phi:=\mathbf 1_{\{\chi_2\}},\qquad \psi:=\mathbf 1_{\{\chi_1,\chi_3\}} ,$$
which are non-negative on $\Sigma$ and have symmetric Fourier coefficients, so that \eqref{Phisym} applies to both. Their expansions in characters are
$\phi=\tfrac14\sum_{t=0}^3(-1)^{t}\xi_t$ and
$\psi=\tfrac12\big(\xi_0-\xi_2\big)$, so
$$T_\phi=\tfrac14\big(I-T+T^2-T^3\big),\qquad T_\psi=\tfrac12\big(I-T^2\big).$$
A direct computation gives
$$T_\phi f(x) = \begin{cases}
 -\frac{1}{4} & x=0,2 \\
 \frac{1}{4} & x=1,3 
\end{cases}, \qquad T_\psi f(x) = \begin{cases}
 \frac{1}{2} & x=0 \\
 -\frac{1}{2} & x=2\\
 0 & x=1,3 
\end{cases}. $$
Hence, by \eqref{Phisym},
$$\Phi(\phi,\psi)=\underset{x\in X}{\mathbb{E}}\ f(x)\,(T_\phi f)(x)\,(T_\psi f)(x)
=\tfrac14\Big(1\cdot(-\tfrac14)\cdot\tfrac12\Big)=-\tfrac{1}{32}<0 . \qedhere$$
\end{proof}
\subsection{Positivity on the cone of non-negative Fourier coefficients}
There is, however, a natural family on which positivity does hold. That is, the cone of functions with non-negative Fourier coefficients.
$$\mathcal{C}:=\Big\{\phi=\sum_{\gamma\in\Gamma}c_\gamma\,\xi_\gamma\ :\ c_\gamma\geq0,\
\textstyle\sum_\gamma c_\gamma<\infty\Big\}\subseteq C(\Sigma).$$ Note that the assumption $\sum_\gamma c_\gamma < \infty$ ensures that the series converges uniformly. Thus, $\phi$ is indeed continuous. Furthermore,
$\phi\in\mathcal{C}$ if and only if $T_\phi$ is a non-negative combination of the
$T_\gamma$. In this case, $T_\phi$ is positivity preserving (i.e., $f\geq 0\Rightarrow T_\phi f \geq 0$).

\subsection{Extending Chu's inequality}
A well known application of H\"older's inequality (see \cite[Lemma 1.6]{Chu}) allows one to bound from below the integral associated with a product of functions that are all the conditional expectations of the same function $f$ under different factors. Here, we will need to generalize this result to a richer class of operators known as doubly-stochastic (see e.g., \cite{Kim}) and more specifically the subclass of those that is positive semi-definite (PSD for short). 

\begin{definition}[PSD doubly-stochastic operator]\label{ds:def}
    A linear operator $M$ on $L^2(\mu)$ is called a \emph{doubly-stochastic operator} if all of the following conditions are satisfied:
    \begin{itemize}
    \item[(1)] $Mf\geq 0$ whenever $f\geq 0$;
    \item[(2)] $M\bm{1}_X=\bm{1}_X$;
    \item[(3)] $M$ is measure-preserving (i.e., $\int_X Mf\,d\mu = \int_X f\,d\mu$).
    \end{itemize}
    It is a \emph{PSD doubly-stochastic operator} if in addition it is self adjoint and 
    \begin{itemize}
        \item[(4)] $\langle Mf,f\rangle \geq 0$ for all $f\in L^2(\mu)$.
    \end{itemize}
\end{definition}
\begin{remark}[Inconsistency of terminology]\label{inconsistent}
Brown \cite{brown} calls operators satisfying $(1)$--$(3)$ \emph{Markov operators}, and shows that they correspond to doubly-stochastic measures. The same class is called doubly-stochastic by Kim \cite{Kim} and by Ryff \cite{Ryff}. We use the latter name because the terminology \emph{Markov operator} is also used in the literature for the strictly larger classes defined by $(1)$ and $(2)$ alone, or by $(1)$ and $(3)$ alone. Note also that for self-adjoint $M$ condition $(3)$ is automatic, since it is equivalent to $M^{*}\bm 1_X=\bm 1_X$.
\end{remark}

Let $\phi\in \mathcal{C}$, non-negative and $\phi(\chi_{\mathrm{triv}})=1$. We show that $M=T_\phi$ is a doubly-stochastic operator. Property $(1)$ follows from the assumption that $\phi\in\mathcal{C}$. Write $\phi = \sum_{\gamma\in \Gamma} c_\gamma \xi_\gamma$, then $\phi(\chi_\mathrm{triv})=1$ implies that $\sum_{\gamma\in \Gamma} c_\gamma = 1$. Hence, $T_{\phi} \bm{1}_X = \sum_{\gamma} c_\gamma T^\gamma \bm{1}_X = \bm{1}_X$, giving $(2)$. Since $\phi$ is real valued, $T_\phi$ is self-adjoint and so
$$\langle T_\phi f,\bm{1}_X\rangle = \langle f,T_\phi \bm{1}_X\rangle  = \int_X f\,d\mu.$$ Finally, $(4)$ follows from the non-negativity of $\phi$.

The main result of this section is inspired by a version of the H\"older's inequality for conditional expectations (see \cite[Lemma 1.6]{Chu}).
\begin{theorem}[Quantitative Positivity]\label{conj}
Let $(X,\mathcal{B},\mu)$ be a probability space, let $k\geq 1$ and let $M_1,\dots,M_k$ be PSD doubly-stochastic operators. Then for all bounded non-negative $f$ we have
$$\int_X f \prod_{i=1}^k M_i f\,d\mu \geq\ \left(\int f\,d\mu\right) ^{k+1}.$$
\end{theorem}
The rest of the section proves Theorem \ref{conj}.
\begin{lemma}[Pointwise Cauchy-Schwarz]\label{pwCS}
Let $(X,\mathcal{B},\mu)$ be a probability space and $M$ a doubly-stochastic operator on $L^2(\mu)$. For every $f,g\in L^\infty(\mu)$ we have
$$|M(fg)|^2 \leq M(f^2)\cdot M(g^2).$$
\end{lemma}
\begin{proof}
    Let $t\in \R$ and observe that $(tf+g)^2\geq 0$. By Property $(1)$ in Definition \ref{ds:def} we conclude that
    $0\leq M((tf+g)^2) = t^2 M(f^2) + 2tM(fg) + M(g^2)$. More precisely, this inequality may only hold almost everywhere. Thus, we shall first take $t\in \Q$ and use the fact that a countable intersection of co-null sets is co-null. Then, discarding the null-set we may pass from $\Q$ to $\R$ via continuity. Viewing this as a quadratic in $t$, the discriminant must be non-positive, giving the desired claim.
\end{proof}
We say that a real function $f$ is \emph{positively $1$-bounded} if it takes values in $[\varepsilon,1]$ for some $0<\varepsilon<1$ (arbitrarily small). Then $\varepsilon\leq  Mf$ from $(1)$, and $0\leq M(\bm{1}_{X}-f) = 1-Mf$ from $(1)$ and $(2)$. Therefore, $\varepsilon\leq Mf\leq 1$ and we may divide by $Mf$ whenever necessary. 
\begin{lemma}\label{PSDused}
    Let $(X,\mathcal{B},\mu)$ be a probability space, let $M$ be a self-adjoint doubly-stochastic operator on $L^2(\mu)$. Then for all positively $1$-bounded $f$ and every $g\in L^2_{\R}(\mu)$, we have 
    $$\langle Mg,g\rangle \leq \int \frac{Mf}{f} g^2\, d\mu.$$
\end{lemma}
\begin{proof}
    Suppose first that $g\in L_{\R}^\infty(\mu)$. Writing $g=\frac{g}{\sqrt{f}}\cdot \sqrt{f}$ the pointwise Cauchy-Schwarz (Lemma~\ref{pwCS}) gives
    $$|Mg| \leq \sqrt{M \left(\frac{g^2}{f}\right)}\cdot \sqrt{M(f)}.$$ Now, by the Cauchy-Schwarz inequality in $L^2(\mu)$ we have
   $$\langle Mg,g\rangle\le\int |g|\sqrt{M \left(\frac{g^2}{f}\right)}\cdot \sqrt{M(f)}\,d\mu
\le\sqrt{\Big(\int g^2\,\tfrac{M(f)}{f}\,d\mu\Big)}\cdot \sqrt{\Big(\int f\,M\left(\tfrac{g^2}{f}\right)\,d\mu\Big)}$$ Since $M$ is self-adjoint, we deduce that $\int f\cdot M\left(\frac{g^2}{f}\right)\,d\mu=\langle \frac{g^2}{f},Mf\rangle=\int g^2\,\frac{Mf}{f}\,d\mu$. In other words, the two components in the inequality above are the same. This removes the square and gives the desired result. Now to extend from $g\in L^\infty_{\R}(\mu)$ to $L^2_{\R}(\mu)$ one may use the fact that the former is dense in the latter and the difference $\int \frac{Mf}{f}g^2\,d\mu - \langle Mg,g\rangle$ is continuous in $g$.
\end{proof}
We now take advantage of the positive semi-definiteness of $M$, which enters through the square root $R(M)$ below. This finally leads to the following inequality, which is the key property used in Chu's proof of \cite[Lemma 1.6]{Chu}.
\begin{lemma}\label{main:lem}
Let $(X,\mathcal{B},\mu)$ be a probability space and let $M$ be a PSD doubly-stochastic operator. Then for every positively $1$-bounded $f$ we have $$\int \frac{f}{Mf}\, d\mu \leq 1.$$
\end{lemma}
\begin{proof}
   Assume first that $M\succeq \eta I$ for some $\eta>0$, so that $M$ is boundedly invertible and so is its positive square root $M^{1/2}$. Let $D$ denote the operator which multiplies by $\frac{Mf}{f}$. Since $f$ is positively $1$-bounded we have $\varepsilon\leq f\leq 1$ and $\varepsilon\leq Mf\leq 1$, hence $\frac{Mf}{f}\in[\varepsilon,\varepsilon^{-1}]$ and $D\succeq\varepsilon I$ is boundedly invertible as well. The previous lemma establishes $\langle Mg,g\rangle \leq \langle Dg,g\rangle$ for all real-valued $g\in L^2(\mu)$, and hence for all complex-valued $g$. Indeed, writing $g=g_1+ig_2$ with $g_1,g_2$ real, both $M$ and $D$ preserve real functions and are self-adjoint, so $\langle Mg,g\rangle=\langle Mg_1,g_1\rangle+\langle Mg_2,g_2\rangle$ and likewise for $D$. We write $M\preceq D$.  Since $D$ is self-adjoint and positive it admits a square root $D^{1/2}$, the same is true for $D^{-1}$. We conclude that $D^{-1/2} M D^{-1/2}\preceq I.$ Set $T=M^{1/2}D^{-1/2}$ and observe that $T^* T = D^{-1/2}M^{1/2}M^{1/2}D^{-1/2} = D^{-1/2} M D^{-1/2}$. Since $T^* T \preceq I$ we have $\|T\|_{\op}\leq 1$, hence $\|T^*\|_{\op}\leq1$ and so we must also have that $TT^* = M^{1/2} D^{-1} M^{1/2}\preceq I.$ Conjugating by $M^{-1/2}$ gives that $D^{-1} \preceq M^{-1}.$ This interchanges the inequality from before. Since $M\bm{1}_X=\bm{1}_X$ we also have $M^{-1}\bm{1}_X=\bm{1}_X$, and therefore
   $$\int \frac{f}{Mf}\,d\mu = \langle D^{-1} \bm{1}_X,\bm{1}_X\rangle \leq \langle M^{-1}\bm{1}_X,\bm{1}_X\rangle = 1.$$

   For a general PSD doubly-stochastic $M$ set $M_\eta:=(1-\eta)M+\eta I$ for $0<\eta<1$. Then $M_\eta$ is again positivity preserving, satisfies $M_\eta\bm 1_X=\bm 1_X$ and $\int M_\eta g\,d\mu=\int g\,d\mu$, is self-adjoint, and $M_\eta\succeq\eta I$. Therefore, the above applies and gives $\int f/(M_\eta f)\,d\mu\leq1$. Moreover, $M_\eta f\geq \varepsilon$ and
   $\|M_\eta f-Mf\|_{L^\infty(\mu)}=\eta\|f-Mf\|_{L^\infty(\mu)}\leq\eta$, so $f/(M_\eta f)\rightarrow f/(Mf)$ uniformly, with all these functions bounded by $\varepsilon^{-1}$. Letting $\eta\rightarrow0$ and using the dominated convergence theorem completes the proof.
\end{proof}
We can now prove Theorem~\ref{conj} using exactly the same argument as in \cite[Lemma 1.6]{Chu}.

\begin{proof}[Proof of Theorem \ref{conj}]
by homogeneity we may assume that $\|f\|_{L^\infty(\mu)}\leq 1$. We may then restrict further to the case that the function $f$ is positively $1$-bounded. Indeed, for the general case, it suffices to
apply the inequality to a scaling of the function $f + \varepsilon$ and take the limit of both
sides when $\varepsilon\rightarrow0$.
Write
$$f = \left( f\cdot \prod_{i=1}^k M_i f\right)^{\frac{1}{k+1}}\cdot \left(\prod_{i=1}^k \frac{f}{M_if}\right)^{\frac{1}{k+1}}.$$ Let $J$ denote the integral on the left hand side of the inequality in Theorem~\ref{conj}. Then, by H\"older inequality we have
$$\left(\int f\, d\mu\right)^{k+1} \leq J \cdot \prod_{i=1}^k \int \frac{f}{M_if}\,d\mu.$$ From the previous lemma, we deduce that the second product is at most $1$, giving the desired inequality.
\end{proof}
In the next section we will use these positivity results to deduce a multiple recurrence theorem in the spirit of Frantzikinakis and Host \cite{FranHost}.
\section{Multiple recurrence for  products of linear terms}\label{prapp}
\begin{definition}
The equation $p(x, y, n) = 0$ is called \emph{partition regular} in $\mathbb{N}$ if for any partition of
$\mathbb{N}$ into finitely many cells, for some $n\in\mathbb{N}$, one of the cells contains distinct $x, y$ that
satisfy the equation.
\end{definition}

In \cite{FranHost}, Frantzikinakis and Host proved the following partition regularity result for quadratic equations.

\begin{theorem} [The three squares theorem]\label{PRFH} Let $p$ be the quadratic form
$$ p(x,y,z) = ax^2 + by^2 + cz^2 + dxy + exz + fyz,$$ where $a, b, c$ are non-zero and $d, e, f$ are arbitrary integers. Suppose that all three forms $$p(x,0,z), p(0,y,z), p(x,x,z)$$ have non-zero square discriminants. Then the equation $p(x,y,n)=0$ is partition regular. 
The last hypothesis means that the three integers
\begin{align*} \nabla_1 &:= e^2 -4ac,\\ \nabla_2&:= f^2 -4bc, \\ \nabla_3&:=(e+f)^2 - 4c(a+b+d)
\end{align*}
are non-zero squares.
\end{theorem}

The main key ingredient in their proof is the following $2$-term multiple recurrence result \cite{FranHostarxiv} (see also \cite[Theorem 2.13]{FranHost} for a generalized result).
\begin{theorem}\label{multrecFH}
Let $l_1$ be positive and $l_2,l_3$ non-negative integers with $l_2\not = l_3$. Then for every set $E\subseteq \mathbb{N}$ of positive multiplicative density, there exist $k_0,m,n\in\mathbb{N}$ such that the integers $k_0\cdot m\cdot (m+l_1n)$ and $k_0\cdot (m+l_2n)\cdot (m+l_3n)$ are distinct
and belong to $E$.
\end{theorem}
In this paper we extend this result to three terms by proving the following $3$-term multiple recurrence theorem.
Following \cite{FranHost}, the pairs $\{l_0,l_1\}$, $\{l_2,l_3\}$ and $\{l_4,l_5\}$, $\{l_6,l_7\}$ play the role of the pairs of shifts in the two products $L_1,L_2$ and $L_1',L_2'$ below. We call the integers $l_0,\dots,l_7$ \emph{admissible} if they are non-negative and
\begin{itemize}
\item[(i)] $l_0,l_1,l_2,l_3$ are pairwise distinct, and
\item[(ii)] $l_4,l_5,l_6,l_7$ are pairwise distinct.
\item[(iii)] $\min\{l_0,l_1,l_2,l_3\} =\min\{l_4,l_5,l_6,l_7\}=0.$
\end{itemize}
\begin{theorem}\label{multrec}
Let $l_0,\dots,l_7$ be admissible. Then for every set $E$ of positive multiplicative density there exist $k_0,m,n,m',n'\in \mathbb{N}$ such that the integers 
\begin{align*}
    k_0\cdot (m+l_0n)\cdot (m+l_1n) \cdot (m'+l_4n')\cdot (m'+l_5n'); \\
   k_0\cdot (m+l_0n) \cdot (m+l_1n) \cdot  (m'+l_6n')\cdot(m'+l_7n');
    \intertext{ and }
     k_0\cdot (m+l_2n)\cdot (m+l_3n)\cdot (m'+l_6n')\cdot (m'+l_7n'),
\end{align*}
are distinct and belong to $E$.
\end{theorem}
Before stating our corollary, we record the parametrization we shall use.

Let $p$ be a quadratic form satisfying the properties of Theorem \ref{PRFH}. By \cite[Proposition 2.11]{FranHost} there are integers $\ell_0,\dots,\ell_4$ with $\ell_0$ positive, $\ell_1\neq\ell_2$, $\ell_3\neq\ell_4$ and $\{\ell_1,\ell_2\}\neq\{\ell_3,\ell_4\}$, such that for all $\kappa,m,n\in\Z$ the integers
$$x:=\kappa\ell_0(m+\ell_1n)(m+\ell_2n),\qquad y:=\kappa\ell_0(m+\ell_3n)(m+\ell_4n)$$
satisfy $p(x,y,\lambda)=0$ for some $\lambda\in\Z$. Writing $\ell_{\min}:=\min\{\ell_1,\ell_2,\ell_3,\ell_4\}$, we call
$$(\ell_1-\ell_{\min},\ \ell_2-\ell_{\min},\ \ell_3-\ell_{\min},\ \ell_4-\ell_{\min})$$
a \emph{shift quadruple} of $p$.

We can now state the \emph{simultaneous partition regularity} counterpart of Theorem~\ref{PRFH}.
\begin{theorem}\label{PR}
Let $p_1,p_2$ be two quadratic forms, each satisfying the properties in Theorem \ref{PRFH}, with shift quadruples $(l_0,l_1,l_2,l_3)$ and $(l_4,l_5,l_6,l_7)$ respectively, and suppose that $l_0,\dots,l_7$ is admissible. Then for any partition of $\mathbb{N}$ into finitely many cells, there exist some $n,n',k\in\mathbb{N}$ and $x,y,x',y'\in \mathbb{N}$ with $x\neq y$ and $x'\neq y'$, so that $p_1(x,y,n)=0$ and $p_2(x',y',n')=0$ and $\frac{x\cdot x'}{k}, \frac{x\cdot y'}{k}, \frac{y\cdot y'}{k}$ are distinct integers which belong to the same cell.
\end{theorem}
\begin{remark}
    For instance, Theorem~\ref{PR} is applicable to the following polynomials  from Frantzikinakis and Host paper.
   $$p_1(x,y,z) = 16x^2+9y^2-z^2. \qquad p_2(x,y,z)=x^2+y^2-xy-z^2.$$ Here, the shift quadruples are $(3,6,4,0)$ and $(1,3,0,2)$, respectively.
\end{remark}
We will now show that Theorem~\ref{multrec} implies Theorem~\ref{PR}.
\begin{proof}
Let $\ell_0,\dots,\ell_4$ be as above for $p_1$, and let $\ell_0',\dots,\ell_4'$ be the corresponding integers for $p_2$, with $x',y'$ defined analogously from $\kappa',m',n'$, so that $p_2(x',y',\lambda')=0$ for some $\lambda'\in\Z$. (The integers $\lambda,\lambda'$ are the ones denoted $n,n'$ in the statement.) Let $(l_0,l_1,l_2,l_3)$ and $(l_4,l_5,l_6,l_7)$ be the shift quadruples of $p_1$ and of $p_2$, so that
\begin{align*}&l_0=\ell_1-\ell_{\min},\quad l_1=\ell_2-\ell_{\min},\quad l_2=\ell_3-\ell_{\min},\quad l_3=\ell_4-\ell_{\min},\\&
l_4=\ell_1'-\ell_{\min}',\quad l_5=\ell_2'-\ell_{\min}',\quad l_6=\ell_3'-\ell_{\min}',\quad l_7=\ell_4'-\ell_{\min}' .
\end{align*}
Replacing $m$ by $m-\ell_{\min}n$ and $m'$ by $m'-\ell_{\min}'n'$ merely relabels the pairs $(m,n)$ and $(m',n')$, and after these substitutions we have $x=\kappa\ell_0L_1(m,n)$, $y=\kappa\ell_0L_2(m,n)$, $x'=\kappa'\ell_0'L_1'(m',n')$ and $y'=\kappa'\ell_0'L_2'(m',n')$ in the notation of Theorem \ref{multrec}.

Let $\mathbb{N}$ be partitioned into finitely many cells. One of them, say $E$, has positive multiplicative density, so Theorem \ref{multrec} provides $k_0,m,n,m',n'\in\mathbb{N}$ for which the three integers $k_0L_1L_1'$, $k_0L_1L_2'$ and $k_0L_2L_2'$ are distinct and belong to $E$. Taking $\kappa:=k_0$, $\kappa':=1$ and $k:=\ell_0\ell_0'$ we obtain
$$\frac{x\cdot x'}{k}=k_0L_1L_1',\qquad \frac{x\cdot y'}{k}=k_0L_1L_2',\qquad \frac{y\cdot y'}{k}=k_0L_2L_2',$$
which are therefore distinct integers lying in the same cell.

Finally, $x\neq y$ and $x'\neq y'$. Indeed, the distinctness assertion of Theorem \ref{multrec} gives $k_0L_1L_2'\neq k_0L_2L_2'$ and $k_0L_1L_1'\neq k_0L_1L_2'$ at the chosen quadruple, whence $L_1\neq L_2$ and $L_1'\neq L_2'$ there. Multiplying by $k_0\ell_0$ and by $\ell_0'$, respectively gives $x\neq y$ and $x'\neq y'$.
\end{proof}
Therefore, in the following sections we focus on proving Theorem \ref{multrec}.
\subsection{Frantzikinakis and Host decomposition of multiplicative functions}
We need some notations. Given a function $f:\mathbb{Z}/N\mathbb{Z}\rightarrow \mathbb{C}$, for some $N\in\mathbb{N}$ and a subset $A\subseteq \mathbb{Z}/N\mathbb{Z}$, we denote the average of $f$ in $A$ by $\mathbb{E}_{n\in A} f(n) = \frac{1}{|A|}\sum_{n\in A} f(n).$ The Gowers uniformity norms of $f$ are defined as follows.
\begin{definition}
Let $f:\mathbb{Z}/N\mathbb{Z}\rightarrow \mathbb{C}$ be a bounded function and let $\mathcal{C}$ denote the complex conjugation. The Gowers $d$-norm\footnote{$\|\cdot\|_{U^d}$ is a seminorm when $d=1$, and a norm for $d\geq 2$.} is defined by the formula
$$\Vert f\Vert _{U^{d}}^{2^{d}}=\mathbb{E} _{x,h_{1},\ldots ,h_{d}\in \mathbb{Z}/N\mathbb{Z}}\prod _{\omega _{1},\ldots ,\omega _{d}\in \{0,1\}}\mathcal{C}^{\omega _{1}+\cdots +\omega _{d}}f\left({x+h_{1}\omega _{1}+\cdots +h_{d}\omega _{d}}\right)\ .$$
\end{definition}

Throughout, $l = \sum_{i=0}^7 l_i$. Given $N\in\mathbb{N}$, we let $\tilde{N}$ denote the smallest prime that is larger than $10l\cdot N$. By Bertrand's postulate, $10lN<\tilde{N}<20lN$. A function $\chi:\mathbb{N}\rightarrow S^1$ is called \emph{multiplicative} if $\chi(n\cdot m) = \chi(n)\cdot \chi(m)$ for all $m,n\in\mathbb{N}$. We write $\mathcal{M}$ for the set of all such functions. For any such function and any $N\in\mathbb{N}$, we denote by $\chi_N:\mathbb{Z}/\tilde{N}\mathbb{Z}\rightarrow\mathbb{C}$ the map defined by 
$$\chi_N(n) = \begin{cases} \chi(n) & 1\leq n\leq N \\ 0 & \text{otherwise} \end{cases}.$$
\begin{definition}
A \emph{kernel} on $\mathbb{Z}/N\mathbb{Z}$ is a non-negative function $\psi:\mathbb{Z}/N\mathbb{Z}\rightarrow \mathbb{R}_{\geq 0}$ with average $1$.
\end{definition}
A key result in our proof is the following decomposition theorem of Frantzikinakis and Host \cite[Theorem 1.6]{FranHostarxiv}.
\begin{theorem}[Structure theorem for multiplicative functions]\label{FHdecom}
Let $\varepsilon>0$, let $\lambda$ be a probability measure on the set of all multiplicative functions $\mathcal{M}$\footnote{Equipped with the pointwise multiplication and the compact-open topology.} and $F:\mathbb{N}\times \mathbb{N}\times \mathbb{R}^+\rightarrow \mathbb{R}^+$ be arbitrary. Then for every sufficiently large $N\in\mathbb{N}$, depending only on $\varepsilon$ and $F$, there exist positive integers $Q$ and $R$, bounded by a constant $K=K(\varepsilon,F)$ which depends only on $\varepsilon$ and $F$ and not on $N$, such that for every $\chi\in \mathcal{M}$ the function $\chi_N = \chi\cdot 1_{[N]}$ admits the decomposition
$$\chi_N (n) = \chi_{N,st}(n) + \chi_{N,un}(n) + \chi_{N,er}(n)$$ for every $n\in \mathbb{Z}/\tilde{N}\mathbb{Z}$, where $\chi_{N,st}, \chi_{N,un}, \chi_{N,er}$ satisfy the following properties:
\begin{itemize}
    \item[(i)] $\chi_{N,st} = \chi_N\ast \psi_{N,1}$ and $\chi_{N,st} + \chi_{N,er} = \chi_N\ast \psi_{N,2}$, where $\psi_{N,1}$, $\psi_{N,2}$ are kernels on $\mathbb{Z}_{\tilde{N}}$ that do not depend on $\chi$, and the convolution product is defined in $\mathbb{Z}/\tilde{N}\mathbb{Z}$. As a consequence, $\chi\mapsto \chi_{N,un}$, $\chi\mapsto \chi_{N,st}$ and $\chi\mapsto\chi_{N,er}$ are continuous, $|\chi_{N,st}|\leq 1$ and $|\chi_{N,un}|, |\chi_{N,er}|\leq 2$;
    \item[(ii)] $|\chi_{N,st}(n+Q)-\chi_{N,st}(n)|\leq \frac{R}{\tilde{N}}$ for every $n\in\mathbb{Z}/\tilde{N}\mathbb{Z}$, where $n+Q$ is taken mod $\tilde{N}$;
    \item[(iii)] $\|\chi_{N,un}\|_{U^{3}(\mathbb{Z}/\tilde{N}\mathbb{Z})}\leq \frac{1}{F(Q,R,\varepsilon)}$;
    \item[(iv)] $$\mathbb{E}_{n\in \mathbb{Z}/\tilde{N}\mathbb{Z}} \int_{\mathcal{M}} |\chi_{N,er}(n)| d\lambda(\chi)\leq \varepsilon.$$
\end{itemize}
\end{theorem}
\begin{remark}\label{pigeonhole}
In \cite{FranHostarxiv} the integers $Q$ and $R$ are allowed to depend on $N$, subject only to the uniform bound $Q,R\leq K(\varepsilon,F)$. Since they are positive integers, only finitely many pairs $(Q,R)$ occur, so there is a single pair $(Q,R)$ and an infinite set $\mathcal{N}\subseteq\mathbb{N}$ of scales for which the conclusion of Theorem \ref{FHdecom} holds with that pair, for every $N\in\mathcal{N}$. In the rest of the paper we work only with scales $N,N'\in\mathcal{N}$, so that $Q$, $R$ and $\eta:=\varepsilon/(QR)$ are genuine constants, depending only on $\varepsilon$ and $F$. 
\end{remark}
\begin{remark}
    A version of this theorem for higher order uniformity norms was also established by Frantzikinakis and Host in \cite{FranHost}. One may see that the very same argument is applicable here. However, performing this argument would require introducing additional parameters and this does not seem to yield much stronger partition regularity results. Therefore, we chose to pursue the modest version introduced above. 
\end{remark}
\subsection{Spectral reformulating of Theorem \ref{multrec}}\label{spectralreform}
In this section we reformulate Theorem \ref{multrec} in the language of multiple correlation sequences. We follow closely the arguments in \cite{FranHost}, and specifically from the preprint which preceded it \cite{FranHostarxiv}, where the main difference is that we apply Theorem~\ref{main:thm} instead of the spectral theorem. 
\begin{definition}[Multiplicative density]
A \emph{multiplicative F\o lner sequence} is an increasing family of finite subsets $(\Phi_N)_{N\in\mathbb{N}}$ of $\mathbb{N}$ satisfying $$\limsup_{N\rightarrow\infty} \frac{|a\cdot \Phi_N \triangle \Phi_N|}{|\Phi_N|} = 0$$ for all $a\in \mathbb{N}$.
The \emph{multiplicative density} $d_{\mathrm{mult}}(E)$ of a subset $E\subseteq \mathbb{N}$, relatively to a F\o lner sequence $\Phi_N$ is defined by the formula
$$d_{\mathrm{mult}}(E):=\limsup_{N\rightarrow\infty} \frac{|E\cap \Phi_N|}{|\Phi_N|}.$$
\end{definition}
Throughout we fix some (any) multiplicative F\o lner sequence $\Phi_N$, and we implicitly assume that $d_{\mathrm{mult}}$ is defined relatively to this sequence.\\

In \cite{FranHost}, Frantzikinakis and Host coined the term \emph{action by dilation}, which in our language is simply a $(\mathbb{Q}^+,\cdot)$-system. Here, $\mathbb{Q}^+$ denote the set of all positive rational numbers. Note that the prime decomposition gives rise to an isomorphism $(\mathbb{Q}^+,\cdot)\cong \bigoplus_{i=1}^\omega\mathbb{Z}$. In particular, we see that the Pontryagin dual $\mathcal{M}$ of $(\mathbb{Q}^+,\cdot)$, which consists of all multiplicative functions on $\mathbb{N}$ and is equipped with pointwise multiplication and the compact open topology, is isomorphic as a topological space to the infinite dimensional torus.\\

The Furstenberg correspondence principle allows us to translate our combinatorial problem into a question about the multiple recurrence of a certain $(\mathbb{Q}^+,\cdot)$-system.
\begin{proposition}[Furstenberg Correspondence Principle]
Let $E$ be a subset of $\mathbb{N}$ and let $\Gamma = (\mathbb{Q}^+,\cdot)$. There exists a $\Gamma$-system $\mathrm{X}=(X,\mathcal{B},\mu,T)$ and a set $A\in\mathcal{B}$ with $\mu(A) = d_{\mathrm{mult}}(E)$, such that for every $r\in \mathbb{N}$ and $n_1,\dots,n_{r}\in \mathbb{N}$, we have
\begin{equation}\label{FCP}
    d_{\mathrm{mult}}(n_1^{-1}E\cap\dots\cap n_r^{-1}E)\geq \mu(T_{n_1}^{-1}A\cap\dots \cap T_{n_{r}}^{-1}A).
\end{equation}
\end{proposition}
By the Furstenberg correspondence principle, we get the following ergodic theoretical reformulation of Theorem \ref{multrec}.
Fix integers $l_0,l_1,\dots,l_7$ as in Theorem \ref{multrec} and let $L_1(m,n) := (m+l_0n)\cdot (m+l_1n)$, $L_2(m,n) := (m+l_2n)\cdot (m+l_3n)$, $L_1'(m,n) := (m+l_4n)\cdot (m+l_5n)$ and $L_2'(m,n) := (m+l_6n)\cdot (m+l_7n)$.
\begin{theorem}[Dynamical reformulation]\label{dynamical}
Let $\mathrm{X}=(X,\mathcal{B},\mu,T)$ be an action by dilation, let $A$ be a measurable set
with $\mu(A)>0$, and set
\begin{align*}&\Theta_N = \{(m,n)\in[N]\times[N] : 1\leq m+l_i n\leq N \text{ for all } 0\leq i \leq 3\},\\&\Theta'_{N'} = \{(m',n')\in[N']\times [N'] : 1\leq m'+l_i n' \leq N' \text{ for all } 4\leq i \leq 7\}.
\end{align*}
Then there exists $\theta>0$ such that, for arbitrarily large $N,N'\in\mathbb{N}$,
$$\underset{(m,n)\in \Theta_N}{\mathbb{E}}~\underset{(m',n')\in \Theta'_{N'}}{\mathbb{E}}\
\mu\big(T^{-1}_{L_1(m,n) L_1'(m',n')} A \cap T^{-1}_{L_1(m,n) L_2'(m',n')}A \cap
T^{-1}_{L_2(m,n) L_2'(m',n')} A\big)\ \geq\ \theta .$$
\end{theorem}
Let us first check that this theorem implies Theorem \ref{multrec}.

We begin by showing that the three products are distinct as polynomials. Since
$l_0,\dots,l_7$ are admissible, $\{l_0,l_1\}\neq\{l_2,l_3\}$ and $\{l_4,l_5\}\neq\{l_6,l_7\}$,
so $L_1\neq L_2$ and $L_1'\neq L_2'$ as elements of $\Z[m,n]$ and $\Z[m',n']$ respectively.
Moreover $L_1L_1'\neq L_2L_2'$, since these polynomials lie in disjoint sets of variables and
are monic in $m^2$ and in $m'^2$, so an identity $L_1L_1'=L_2L_2'$ would force $L_1=L_2$.
Consequently the three sets
$$\{L_1'=L_2'\},\qquad \{L_1=L_2\},\qquad \{L_1L_1'=L_2L_2'\}$$
are zero sets of non-zero polynomials. A non-zero polynomial in two variables vanishes at
$O(N)$ points of $[N]^2$, and one in four variables at $O(NN'^2+N^2N')$ points of
$[N]^2\times[N']^2$. Since $|\Theta_N|\geq cN^2$ and $|\Theta'_{N'}|\geq c' N'^2$ for some constants $c,c'$ depending only on $l_0,\dots,l_7$, the union of the three sets meets $\Theta_N\times\Theta'_{N'}$ in a subset of
density $O(1/N)+O(1/N')$.

Let $\theta>0$ be as in the Dynamical reformulation, and choose $N,N'$ large enough that
this density is smaller than $\theta$ and that the conclusion of that theorem holds for
$N,N'$. The set of quadruples in $\Theta_N\times\Theta'_{N'}$ at which
$$\mu\big(T^{-1}_{L_1L_1'} A \cap T^{-1}_{L_1L_2'}A \cap T^{-1}_{L_2L_2'} A\big)>0$$
has density at least $\theta$ there, so we may fix such a quadruple $(m,n,m',n')$ lying
outside all three sets above. Apply  \eqref{FCP} with $r=3$ and
$$n_1=L_1(m,n)L_1'(m',n'),\qquad n_2=L_1(m,n)L_2'(m',n'),\qquad n_3=L_2(m,n)L_2'(m',n').$$
The left-hand side of \eqref{FCP} is positive by the choice of the quadruple, so the set
$n_1^{-1}E\cap n_2^{-1}E\cap n_3^{-1}E$ is non-empty. Any $k_0$ in it satisfies
$k_0n_1,k_0n_2,k_0n_3\in E$, and these are exactly the three integers displayed in
Theorem \ref{multrec}. Finally $n_1\neq n_2$, $n_2\neq n_3$ and $n_1\neq n_3$, since the
quadruple avoids $\{L_1'=L_2'\}$, $\{L_1=L_2\}$ and $\{L_1L_1'=L_2L_2'\}$ respectively, and
$L_1,L_2'$ are positive there. We conclude that $k_0n_1,k_0n_2,k_0n_3$ are distinct, as required.

It is left to prove this theorem.\\

Let $R(m,n) := \frac{L_2(m,n)}{L_1(m,n)}$ and $R'(m,n) := \frac{L_2'(m,n)}{L_1'(m,n)}$. Since $T$ is measure preserving, $T^{-1}_{L_1L_1'}A\cap T^{-1}_{L_1L_2'}A\cap T^{-1}_{L_2L_2'}A$
and $A\cap T^{-1}_{R'(m',n')}A\cap T^{-1}_{R(m,n)R'(m',n')}A$ have the same measure, so
Theorem \ref{dynamical} is equivalent to the existence of $\theta>0$ such that, for arbitrarily
large $N,N'$,
$$\underset{(m,n)\in\Theta_N}{\mathbb{E}}~\underset{(m',n')\in\Theta'_{N'}}{\mathbb{E}}\
\mu\big(A\cap T_{R'(m',n')}^{-1} A \cap T_{R(m,n)\cdot R'(m',n')}^{-1}A\big)\ \geq\ \theta .$$ From \eqref{spec2} (applied with $f=g=h=1_A$ and $S=T$), there exist a Borel probability measure $\lambda=\lambda_A$ on $\mathcal{M}$ and an operator $G=G_A:L^2(\mathcal{M},\lambda_A)\rightarrow L^2(\mathcal{M},\lambda_A)$ so that
$$\mu(A\cap T_{R'(m',n')}^{-1} A \cap T_{R(m,n)\cdot R'(m',n')}^{-1}A) = \int_{\mathcal{M}} G(\xi_{R'(m',n')})(\chi)\cdot \xi_{R(m,n)}(\chi) d\lambda(\chi),$$ where for every $t\in \mathbb{N}$, $\xi_t(\chi) := \chi(t)$ is the evaluation map. Therefore, it suffices to prove that 
\begin{equation}\label{liminf}
\underset{m,n\in \Theta_N}{\mathbb{E}} ~\underset{m',n'\in \Theta'_{N'}}{\mathbb{E}}
\int_{\mathcal{M}} G(\xi_{R'(m',n')})(\chi)\cdot \xi_{R(m,n)}(\chi)~d\lambda(\chi)\ \geq\ \theta
\quad \text{for arbitrarily large } N,N'\in\mathbb{N},
\end{equation}
for some $\theta>0$. That is, for every $N_0$ there are $N,N'>N_0$ for which the displayed average is at least
$\theta$. Indeed, an average $\geq \theta$ produces a set of quadruples $(m,n,m',n')\in\Theta_N\times\Theta'_{N'}$ of density $\geq \theta$ for which the integrand is positive,
where $\Theta_N = \{(m,n)\in[N]\times[N] : 1\leq m+l_i n\leq N \text{ for all } 0\leq i \leq 3\}$, and $\Theta'_{N'} = \{(m',n')\in[N']\times [N'] : 1\leq m'+l_i n' \leq N' \text{ for all } 4\leq i \leq 7\}$.\\

We define the truncated $\Xi_{N,m,n}$ and $\Xi'_{N',m',n'}$ by
\begin{align*}
    \Xi_{N,m,n}(\chi) &= \overline{\chi_N(m+l_0n)}\cdot \overline{\chi_N(m+l_1n)}\cdot \chi_N(m+l_2n)
\cdot \chi_N(m+l_3n),\\
 \Xi'_{N',m',n'}(\chi) &= \overline{\chi_{N'}(m'+l_4n')}\cdot \overline{\chi_{N'}(m'+l_5n')}\cdot \chi_{N'}(m'+l_6n')
\cdot \chi_{N'}(m'+l_7n').
\end{align*}
\begin{lemma}\label{liminf'}
Let $\theta>0$. To prove \eqref{liminf} for that $\theta$, it suffices to show that for arbitrarily large $N,N'\in\mathbb{N}$
 $$\underset{m,n\in \mathbb{Z}/\tilde{N}\mathbb{Z}}{\mathbb{E}} ~\underset{m',n'\in \mathbb{Z}/\tilde{N}'\mathbb{Z}}{\mathbb{E}} \int_{\mathcal{M}}  G(1_{[N']}(n')\Xi'_{N',m',n'})(\chi)\cdot  1_{[N]}(n) \Xi_{N,m,n}(\chi)~d\lambda(\chi)\geq\theta.$$
\end{lemma}
\begin{proof}
First, observe that
\begin{align*}
&\underset{m,n\in \Theta_N}{\mathbb{E}} ~\underset{m',n'\in \Theta'_{N'}}{\mathbb{E}} \int_{\mathcal{M}} G(\xi_{R'(m',n')})(\chi)\cdot \xi_{R(m,n)}(\chi)~d\lambda(\chi) = \\
&\frac{\tilde{N}^2\cdot \tilde{N'}^2}{|\Theta_N|\cdot|\Theta'_{N'}|}\cdot
    \underset{m,n\in \mathbb{Z}/\tilde{N}\mathbb{Z}}{\mathbb{E}} ~\underset{m',n'\in \mathbb{Z}/\tilde{N}'\mathbb{Z}}{\mathbb{E}} \int_{\mathcal{M}}  G(1_{[N']}(n')\Xi'_{N',m',n'})(\chi)\cdot  1_{[N]}(n) \Xi_{N,m,n}(\chi)~d\lambda(\chi).
\end{align*}
For $(m,n)\in\Theta_N$, each of $m+l_0n,\,m+l_1n,\,m+l_2n,\,m+l_3n$ lies in $[N]$, so $\chi_N$ agrees with $\chi$ at each of them, and by complete multiplicativity $\Xi_{N,m,n}(\chi)=\chi(R(m,n))=\xi_{R(m,n)}(\chi)$. Similarly $\Xi'_{N',m',n'}(\chi)=\xi_{R'(m',n')}(\chi)$ for $(m',n')\in\Theta'_{N'}$. On the other hand, if $n>N$, then the term inside the average is zero because $1_{[N]}(n)=0$. By $(iii)$ some $l_{i_0}$ with $0\leq i_0\leq 3$ equals $0$, so the factor
$\chi_N(m+l_{i_0}n)=\chi_N(m)$ forces $m\in[N]$. This, together with $n\in[N]$ and
$m+l_in\in[N]$ says exactly that $(m,n)\in\Theta_N$. Using the transpose to move $G$ to the other side, the same argument shows that the summand vanishes unless $(m',n')\in\Theta'_{N'}$. Hence, the summand vanishes off $\Theta_N\times\Theta'_{N'}$, and the two averages differ by the stated factor. Since $\tilde{N}>10lN$ and $|\Theta_N|\leq N^2$, and similarly for $N'$, the factor
$\tilde{N}^2\tilde{N}'^2/(|\Theta_N||\Theta'_{N'}|)$ is at least $10^4l^4\geq 1$. Hence the
average in \eqref{liminf} is at least the average in the statement, and any lower bound for the
latter transfers to the former.
\end{proof}
The following estimate was established in \cite[Lemma 2.7]{FranHostarxiv}.
\begin{lemma}[$U^3$-estimate]\label{U^3estimate}
    Let $a_i$, $i=0,1,2,3$, be functions on $\mathbb{Z}/\tilde{N}\mathbb{Z}$, with $\|a_i\|_{L^\infty(\mathbb{Z}/\tilde{N}\mathbb{Z})} \leq 1$ and $l_1,l_2,l_3\in\mathbb{N}$ be distinct. Then there exists a constant $c_2$ depending only on $l=l_1+l_2+l_3$ such that $$\left|\underset{m,n\in \mathbb{Z}/\tilde{N}\mathbb{Z}}{\mathbb{E}} 1_{[N]}(n)\cdot a_0(m)\cdot a_1(m+l_1n)\cdot a_2(m+l_2n)\cdot a_3(m+l_3n)\right| \leq c_2 \min_{0\leq j \leq 3} \|a_j\|_{U^3(\mathbb{Z}/\tilde{N}\mathbb{Z})}^{\frac{1}{2}} + \frac{2}{\tilde{N}}.$$
\end{lemma}
By admissibility each quadruple consists of $0$ together with three distinct positive integers. Hence, we may apply this lemma to the four shifts $l_0,l_1,l_2,l_3$ of $\Xi_{N,m,n}$, and to the four shifts $l_4,l_5,l_6,l_7$ of $\Xi'_{N',m',n'}$, after relabeling the four functions so that the one whose shift is $0$ is $a_0$.

For non-negative $\psi$ on $\Z/\tilde{N}\Z$ and $n\in \Z/\tilde{N}\Z$, let $\xi_{\psi,N,n}$
denote the map on $\mathcal{M}$ sending $\chi$ to $\chi_N\ast\psi(n)$, and define
$\xi_{\psi',N',n'}$ similarly. Set
\begin{align*}
\Xi_{\psi,N,m,n} &:= \overline{\xi_{\psi,N,m+l_0n}}\cdot \overline{\xi_{\psi,N,m+l_1n}}\cdot
   \xi_{\psi,N,m+l_2n}\cdot \xi_{\psi,N,m+l_3n},\\
\Xi'_{\psi',N',m',n'} &:= \overline{\xi_{\psi',N',m'+l_4n'}}\cdot \overline{\xi_{\psi',N',m'+l_5n'}}\cdot
   \xi_{\psi',N',m'+l_6n'}\cdot \xi_{\psi',N',m'+l_7n'} .
\end{align*}
The following lemma is inspired by \cite[Lemma 2.8]{FranHostarxiv}.
\begin{lemma}\label{non-negative}
Let $\mathrm{X}=(X,\mathcal{B},\mu,T)$ be an action by dilations, let $A\subseteq X$ be
measurable, and let $G$ and $\lambda$ be as above. Let $\psi,\psi'$ be non-negative
functions on $\Z/\tilde{N}\Z$ and $\Z/\tilde{N}'\Z$ respectively. Then
$$\int_{\mathcal{M}} G\big(\Xi'_{\psi',N',m',n'}\big)(\chi)\cdot \Xi_{\psi,N,m,n}(\chi)~d\lambda(\chi)\ \geq\ 0$$
for all $m,n\in \Z/\tilde{N}\Z$ and $m',n'\in \Z/\tilde{N}'\Z$.
\end{lemma}

\begin{proof}
By definition $\chi_N\ast\psi(n)=\E_{k\in\Z/\tilde{N}\Z}\psi(n-k)\chi_N(k)$, so with
$a_n(k):=\psi(n-k)\geq 0$ we have $\xi_{\psi,N,n}=\E_{k}a_n(k)\,\xi_{k}$, the sum being
supported on $k\in[N]$ since $\chi_N$ vanishes off $[N]$. Writing $a'_{n'}(k'):=\psi'(n'-k')$
similarly, and expanding the eight factors, the left-hand side equals
$$\E_{k_1,\dots,k_4}\ \E_{k'_1,\dots,k'_4}\ \Big(\prod_{i=1}^4 a_{m_i}(k_i)\,a'_{m'_i}(k'_i)\Big)
\int_{\mathcal{M}} G\big(\xi_{t'}\big)(\chi)\cdot \xi_{t}(\chi)~d\lambda(\chi),$$
where $m_1=m+l_0n$, $m_2=m+l_1n$, $m_3=m+l_2n$, $m_4=m+l_3n$, and $m'_1,\dots,m'_4$ are the
corresponding forms in $m',n'$, and where
$$t:=\frac{k_3k_4}{k_1k_2}\in\Q^+,\qquad t':=\frac{k'_3k'_4}{k'_1k'_2}\in\Q^+ .$$
Here we used that $\chi$ is completely multiplicative, so that
$\overline{\xi_{k_1}}\,\overline{\xi_{k_2}}\,\xi_{k_3}\xi_{k_4}=\xi_{t}$.
By \eqref{spec2gen} with $f=g=h=1_A$,
$$\int_{\mathcal{M}} G(\xi_{t'})(\chi)\cdot \xi_{t}(\chi)~d\lambda(\chi)
= \int_X T_{t'}\big(T_{t}1_A\cdot 1_A\big)\cdot 1_A~d\mu
= \mu\big(A\cap T_{t'}^{-1}A\cap T_{t\cdot t'}^{-1}A\big)\ \geq 0 .$$
Since every coefficient $a_{m_i}(k_i)a'_{m'_i}(k'_i)$ is non-negative, the whole expression
is non-negative.
\end{proof}
\subsection{Completing the proof of Theorem \ref{multrec}}\label{completing}
We proved in the previous section that in order to prove Theorem \ref{multrec}, it suffices to show that the term appearing in Lemma \ref{liminf'} is positive for arbitrarily large $N,N'$. The main novelty here (compared to \cite{FranHostarxiv,FranHost}) is the estimate \eqref{l1estimate}, which is where the operator $G$ enters. We let $f=1_A$ denote the characteristic function of $A$, and set  $\delta:=\mu(A) = \int f d\mu$. Let $\varepsilon = c_3\cdot \delta^4$ and $F(x,y,z) = c_4^2\frac{x^4y^4}{z}$ where $c_3,c_4$ are constants depending only on $l_0,\dots,l_7$ and on $\delta$, to be chosen later. Let $$A(N,N'):=\int \left(\underset{m',n'\in\mathbb{Z}/\tilde{N}'\mathbb{Z}}{\mathbb{E}} G\left(1_{[N']}(n')\cdot  \Xi'_{N',m',n'}\right)\right) \cdot  \left(\underset{m,n\in \mathbb{Z}/\tilde{N}\mathbb{Z}}{\mathbb{E}} 1_{[N]}(n)\cdot  \Xi_{N,m,n}\right) d\lambda.$$
We apply Theorem \ref{FHdecom} to the $\varepsilon, F, \lambda$ defined above. Let $Q,R$ be as in the theorem, and write $\xi_{N,n} = \xi_{N,n}^u + \xi_{N,n}^s + \xi_{N,n}^e$ where $\xi_{N,n}^u(\chi) = \chi_{N,un}(n)$, $\xi_{N,n}^s(\chi)=\chi_{N,st}(n)$ and $\xi_{N,n}^e(\chi) = \chi_{N,er}(n)$ satisfy the properties of the theorem. We also write $\xi_{N,n}^{s,e} = \xi_{N,n}^s+\xi_{N,n}^e$. Accordingly, for $\ast\in\{u,s,e,(s,e)\}$ we write $\Xi^{\ast}_{N,m,n}$ and
$\Xi'^{\ast}_{N',m',n'}$ for the functions obtained by replacing each factor
$\xi_{N,\cdot}$ in the definition of $\Xi_{N,m,n}$, respectively
$\Xi'_{N',m',n'}$, by $\xi^{\ast}_{N,\cdot}$. That is,
$$\Xi^{s}_{N,m,n} = \overline{\xi^{s}_{N,m+l_0n}}\cdot \overline{\xi^{s}_{N,m+l_1n}}\cdot
\xi^{s}_{N,m+l_2n}\cdot \xi^{s}_{N,m+l_3n},$$
and similarly for the other options. Now, look at
$$A_1(N,N') :=\int \left(\underset{m',n'\in \mathbb{Z}/\tilde{N}'\mathbb{Z}}{\mathbb{E}} G\left(1_{[N']}(n')\cdot  \Xi'^{s,e}_{N',m',n'}\right)\right) \cdot\left(\underset{m,n\in \mathbb{Z}/\tilde{N}\mathbb{Z}}{\mathbb{E}} 1_{[N]}(n)\cdot  \Xi^{s,e}_{N,m,n}\right) d\lambda.$$
Namely, the term obtained by replacing each instance of $\xi$ with $\xi^{s,e}$ (removing the uniform part). We bound $A(N,N')-A_1(N,N')$. To do so we introduce an intermediate term $$B(N,N'):=\int \left(\underset{m',n'\in \mathbb{Z}/\tilde{N}'\mathbb{Z}}{\mathbb{E}} G\left(1_{[N']}(n')\cdot  \Xi'^{s,e}_{N',m',n'}\right)\right) \cdot\left(\underset{m,n\in \mathbb{Z}/\tilde{N}\mathbb{Z}}{\mathbb{E}} 1_{[N]}(n)\cdot  \Xi_{N,m,n}\right) d\lambda.$$
By Cauchy--Schwarz and since $\lambda$ is a probability measure, we have for all bounded $\phi,\psi$ that $|\int G(\phi)\cdot \psi|\leq \|G\phi\|_{L^2(\lambda)} \|\psi\|_{L^2(\lambda)} \leq \|G\|\cdot \|\phi\|_\infty \cdot \|\psi\|_\infty$. From the bounds on the characters in Theorem~\ref{FHdecom}(i), we deduce that
$$|B(N,N')-A_1(N,N')| \leq \|G\|_{op}\cdot \big\|\underset{m,n\in \mathbb{Z}/\tilde{N}\mathbb{Z}}{\mathbb{E}} 1_{[N]}(n)\cdot \left(\Xi_{N,m,n}-\Xi^{s,e}_{N,m,n}\right)\big\|_{\infty}.$$
Recall that for every $\chi$, $\Xi_{N,m,n}(\chi) = \overline{\chi_N(m+l_0n)}\cdot \overline{\chi_N(m+l_1n)}\cdot \chi_N(m+l_2n)\cdot \chi_N(m+l_3n).$ Therefore, the average on the right hand side in the equation above can be written as a sum of $4$ terms, each is a multiple of $4$ terms, taking the same form as in Lemma \ref{U^3estimate}. Moreover, each of these summands contains at least one multiple that has $U^3$ norm $\leq \frac{1}{F(Q,R,\varepsilon)}$. Since $|\chi_{N,un}|\leq2$ by Theorem \ref{FHdecom}(i), Lemma \ref{U^3estimate} is applied to $\tfrac12\xi^u_{N,\cdot}$, which costs a factor $\sqrt{2}$ that we absorb into $c_2$, and doubles the error term $\frac{2}{\tilde{N}}$ of that lemma to $\frac{4}{\tilde{N}}$. Therefore, we deduce that
$$|B(N,N')-A_1(N,N')| \leq \|G\|_{op} \cdot \frac{4c_2}{F(Q,R,\varepsilon)^{\frac{1}{2}}} + \frac{16}{\tilde{N}}. $$
Using the transpose (i.e., $\int G\phi\cdot \psi\,d\lambda = \int \phi\cdot G^*\psi\,d\lambda$.), we can also write 
$$B(N,N') = \int \left(\underset{m',n'\in \mathbb{Z}/\tilde{N}'\mathbb{Z}}{\mathbb{E}} 1_{[N']}(n')\cdot  \Xi'^{s,e}_{N',m',n'}\right) \cdot\left(\underset{m,n\in \mathbb{Z}/\tilde{N}\mathbb{Z}}{\mathbb{E}} G^{*}\left(1_{[N]}(n)\cdot  \Xi_{N,m,n}\right)\right)\, d\lambda$$
and obtain the estimate
$|A(N,N')-B(N,N')| \leq \frac{4c_2\cdot \|G\|_{op}}{F(Q,R,\varepsilon)^{\frac{1}{2}}} + \frac{16}{\tilde{N'}}$ using the exact same argument as above (and using the well known fact that $\|G\|_{op} = \|G^{*}\|_{op}$). By the triangle inequality we deduce that
\begin{equation}\label{A-A1}|A(N,N')-A_1(N,N')| < \frac{8c_2\cdot \|G\|_{op}}{F(Q,R,\varepsilon)^{\frac{1}{2}}} + \frac{16}{\tilde{N'}} + \frac{16}{\tilde{N}}.\end{equation}
We now work with $A_1(N,N')$. We want to eliminate the error term, but first, we take advantage of the (almost) periodicity of the structured component. Recall from Theorem \ref{FHdecom} that $\chi_{N,st}$ is periodic with period $Q$ (up to the
stated error). Accordingly, for $m,n\in\Z/\tilde{N}\Z$ and $m',n'\in\Z/\tilde{N}'\Z$ we set
\begin{align*}
\tilde{\Xi}_{N,m,n} &:= \Xi_{N,m,Qn}
 = \overline{\xi_{N,m+l_0Qn}}\cdot \overline{\xi_{N,m+l_1Qn}}\cdot \xi_{N,m+l_2Qn}\cdot \xi_{N,m+l_3Qn},\\
\tilde{\Xi}'_{N',m',n'} &:= \Xi'_{N',m',Qn'}
 = \overline{\xi_{N',m'+l_4Qn'}}\cdot \overline{\xi_{N',m'+l_5Qn'}}\cdot
   \xi_{N',m'+l_6Qn'}\cdot \xi_{N',m'+l_7Qn'},
\end{align*}
and we write $\tilde{\Xi}^{\ast}_{N,m,n}$, $\tilde{\Xi}'^{\ast}_{N',m',n'}$ for the
corresponding decorated versions, $\ast\in\{u,s,e,(s,e)\}$, obtained by replacing each
factor $\xi_{N,\cdot}$ by $\xi^{\ast}_{N,\cdot}$ as above. Let $\eta:=\frac{\varepsilon}{QR}$, so that $Q\lfloor\eta N\rfloor\leq N$. Every
summand on the left-hand side below is non-negative by Lemma \ref{non-negative}, and the
summands indexed by $m\in\Z/\tilde{N}\Z$, $1\leq n\leq\lfloor\eta N\rfloor$,
$m'\in\Z/\tilde{N}'\Z$, $1\leq n'\leq\lfloor\eta N'\rfloor$ after the substitutions
$n\mapsto Qn$, $n'\mapsto Qn'$ form a sub-family of them, on which
$1_{[N]}(Qn)=1_{[N']}(Qn')=1$. Discarding the remaining summands therefore gives
\begin{align*}
    &\sum_{m,n\in \mathbb{Z}/\tilde{N}\mathbb{Z}} \sum_{m',n'\in \mathbb{Z}/\tilde{N'}\mathbb{Z}} \int_{\mathcal{M}} G\big(1_{[N']}(n')\,\Xi'^{s,e}_{N',m',n'}\big)(\chi)\cdot 1_{[N]}(n)\,\Xi^{s,e}_{N,m,n}(\chi)~d\lambda(\chi) \\
    &\qquad\geq\ \sum_{m\in \mathbb{Z}/\tilde{N}\mathbb{Z}} \sum_{n=1}^{\lfloor\eta N\rfloor} \sum_{m'\in \mathbb{Z}/\tilde{N'}\mathbb{Z}} \sum_{n'=1}^{\lfloor\eta N'\rfloor} \int_{\mathcal{M}} G\big( \tilde{\Xi}'^{s,e}_{N',m',n'}\big)(\chi)\cdot \tilde{\Xi}^{s,e}_{N,m,n}(\chi)~d\lambda(\chi).
\end{align*}

Indeed, the summands associated with $n>N$ or $n'>N'$ in the first term are zero. For the rest of the terms we notice that the right-hand side has fewer summands and so the inequality follows from Lemma \ref{non-negative}. We denote
$$A_2(N,N') := \underset{m\in \mathbb{Z}/\tilde{N}\mathbb{Z}}{\mathbb{E}} ~\underset{n\leq \lfloor\eta N\rfloor}{\mathbb{E}}~\underset{m'\in \mathbb{Z}/\tilde{N'}\mathbb{Z}}{\mathbb{E}}~\underset{n'\leq \lfloor\eta N'\rfloor}{\mathbb{E}} \int_{\mathcal{M}} G(\tilde{\Xi}'^{s,e}_{N',m',n'})\cdot\tilde{\Xi}^{s,e}_{N,m,n}d\lambda.$$

By the inequality above we have
\begin{equation}\label{A1 geq A2}
\begin{split}
    A_1(N,N')&\geq \frac{\lfloor \eta N\rfloor}{\tilde{N}}\frac{\lfloor \eta N'\rfloor}{\tilde{N'}}\cdot A_2(N,N')\\ &\geq \frac{\eta^2}{500\cdot l^2}\cdot  A_2(N,N')\\ & = \frac{\varepsilon^2}{500\cdot l^2\cdot Q^2\cdot R^2 }\cdot A_2(N,N').
    \end{split}
\end{equation}
Here we used $\tilde{N}<20lN$ and $\tilde{N}'<20lN'$, so that
$\frac{\lfloor\eta N\rfloor}{\tilde N}\geq\frac{\eta N-1}{20lN}$, and similarly for $N'$. The resulting constant is at least $\frac{\eta^2}{500l^2}$ once $N,N'$ are large enough.
We therefore work with $A_2(N,N')$ from now on. We are set to remove the error term. Set $$A_3(N,N') := \underset{m\in \mathbb{Z}/\tilde{N}\mathbb{Z}}{\mathbb{E}} ~\underset{n\leq \lfloor\eta N\rfloor}{\mathbb{E}}~\underset{m'\in \mathbb{Z}/\tilde{N'}\mathbb{Z}}{\mathbb{E}}~\underset{n'\leq \lfloor\eta N'\rfloor}{\mathbb{E}} \int_{\mathcal{M}} G(\tilde{\Xi}'^{s}_{N',m',n'})\cdot\tilde{\Xi}^{s}_{N,m,n}d\lambda.$$
To estimate $|A_2(N,N')-A_3(N,N')|$ we use a similar argument as we used to get \eqref{A-A1}, but here we have to rely on another general estimate involving the $L^1$ norm. By the Cauchy--Schwarz inequality we have  $$\left|\int_{\mathcal{M}} G(\phi) \psi d\lambda \right| \leq \|G\|_{op}\cdot  \|\phi\|_{L^2(\lambda)} \|\psi\|_{L^2(\lambda)}.$$ Since the inequality $\|\cdot \|_{L^2} \leq \sqrt{ \|\cdot \|_{L^1} \cdot \|\cdot \|_\infty}$ holds in all probability spaces we deduce that
\begin{equation}\label{l1estimate}
\left|\int_{\mathcal{M}} G(\phi) \psi d\lambda \right| \leq \|G\|_{op}\cdot \sqrt{\|\phi\|_{L^1(\lambda)}\cdot \|\phi\|_{L^\infty(\lambda)} \cdot \|\psi\|_{L^1(\lambda)}\cdot \|\psi\|_{L^\infty(\lambda)}}.
\end{equation}
Recall that by Theorem \ref{FHdecom}(iv), $\underset{n\in\mathbb{Z}/\tilde{N}\mathbb{Z}}{\mathbb{E}}\|\xi^e_{N,n}\|_{L^1(\lambda)}\leq\varepsilon$. Again we introduce an intermediate term $$B_2(N,N'):= \underset{m\in \mathbb{Z}/\tilde{N}\mathbb{Z}}{\mathbb{E}} ~\underset{n\leq \lfloor\eta N\rfloor}{\mathbb{E}}~\underset{m'\in \mathbb{Z}/\tilde{N'}\mathbb{Z}}{\mathbb{E}}~\underset{n'\leq \lfloor\eta N'\rfloor}{\mathbb{E}} \int_{\mathcal{M}} G( \tilde{\Xi}'^{s,e}_{N',m',n'})\cdot \tilde{\Xi}^{s}_{N,m,n}d\lambda.$$
Since all the $\xi$'s are $1$-bounded, \eqref{l1estimate} implies that
$$|A_2(N,N') - B_2(N,N')| \leq  \|G\|_{op} \underset{m\in \mathbb{Z}/\tilde{N}\mathbb{Z}}{\mathbb{E}} ~\underset{n\leq \lfloor\eta N\rfloor}{\mathbb{E}}\sqrt{\|\tilde{\Xi}^{s,e}_{N,m,n} - \tilde{\Xi}^s_{N,m,n}\|_{L^1}}.$$

Once again, we can write $\tilde{\Xi}^{s,e}_{N,m,n}-\tilde{\Xi}^s_{N,m,n}$ as four summands, each a multiple of $4$ terms, where all terms are $1$-bounded in $L^\infty$ norm, and exactly one of them is of the form $\xi^e_{N,m+l_iQn}$ for some $0\leq i\leq 3$. For each fixed $n$ the map $m\mapsto m+l_iQn$ is a bijection of $\mathbb{Z}/\tilde{N}\mathbb{Z}$, and therefore
$$\underset{m\in \mathbb{Z}/\tilde{N}\mathbb{Z}}{\mathbb{E}}\ \big\|\xi^e_{N,m+l_iQn}\big\|_{L^1(\lambda)} = \underset{m\in \mathbb{Z}/\tilde{N}\mathbb{Z}}{\mathbb{E}}\ \big\|\xi^e_{N,m}\big\|_{L^1(\lambda)}\ \leq\ \varepsilon$$
by Theorem \ref{FHdecom}(iv). By Jensen's inequality ( $\mathbb{E}_n\sqrt{x_n}\leq\sqrt{\mathbb{E}_n x_n}$) and the  Cauchy--Schwarz inequality, this gives the estimate
$$|A_2(N,N') - B_2(N,N')| \leq 4\cdot \sqrt{\varepsilon}\cdot \|G\|_{op}.$$

Using the transpose to move $G$ to the other term, as in the previous argument, we also get the bound
$$|B_2(N,N') - A_3(N,N')|\leq 4 \sqrt{\varepsilon}\cdot \|G^{*}\|_{op},$$ and so by the triangle inequality we have
\begin{equation}\label{A2 - A3}
|A_2(N,N')-A_3(N,N')| < 8 \sqrt{\varepsilon} \cdot\|G\|_{op}.
\end{equation}
It is left to estimate $A_3(N,N')$. Now that we are left with the structure term we can use the periodicity. Recall that $$\tilde{\Xi}^s_{N,m,n} (\chi) = \overline{\chi^s_N(m+l_0Qn)}\cdot \overline{\chi^s_N(m+l_1Qn)}\cdot \chi^s_{N}(m+l_2Qn) \cdot \chi^s_{N}(m+l_3Qn).$$ 
By Theorem \ref{FHdecom}(ii), iterating $l_i n$ times in steps of $Q$ gives
$\|\xi^s_{N,m+l_iQn}-\xi^s_{N,m}\|_\infty \leq l_i\, n\, \frac{R}{\tilde{N}}$.
Since all four factors are $1$-bounded by Theorem \ref{FHdecom}(i), telescoping
yields, for every $m$ and every $1\leq n\leq \lfloor\eta N\rfloor$,
$$\big\|\tilde{\Xi}^s_{N,m,n} - |\xi^s_{N,m}|^4\big\|_{\infty}
\leq l\cdot n\cdot \frac{R}{\tilde{N}}
\leq l\cdot \eta N\cdot \frac{R}{\tilde{N}}
\leq \frac{\varepsilon}{Q},$$
where we used $\eta=\frac{\varepsilon}{QR}$ and $\tilde{N}>l\cdot N$. Similarly
$$\big\|\tilde{\Xi}'^s_{N',m',n'} - |\xi^s_{N',m'}|^4\big\|_{\infty}\leq \frac{\varepsilon}{Q}.$$

Let $$A_4(N,N') := \int_{\mathcal{M}} G(\mathbb{E}_{m'\in \mathbb{Z}/\tilde{N'}\mathbb{Z}} |\xi^s_{N',m'}|^4) \cdot \mathbb{E}_{m\in \mathbb{Z}/\tilde{N}\mathbb{Z}} |\xi^s_{N,m}|^4 d\lambda,$$ by the Cauchy--Schwarz inequality we have
\begin{equation}\label{A3-A4}
    |A_3(N,N')-A_4(N,N')|\leq 16\cdot \|G\|_{op}\cdot \frac{\varepsilon}{Q}  \leq 16\cdot \varepsilon\cdot \|G\|_{op}.
\end{equation}
It is left to bound $A_4$ from below. Here we take advantage of the main positivity result from the previous section. This is also the main step where our proof deviates from \cite{FranHostarxiv,FranHost}.\\
\noindent Write
$$\Psi:=\underset{m\in\mathbb{Z}/\tilde{N}\mathbb{Z}}{\mathbb{E}}\big|\xi^{s}_{N,m}\big|^{4},
\qquad
\Psi':=\underset{m'\in\mathbb{Z}/\tilde{N}'\mathbb{Z}}{\mathbb{E}}\big|\xi^{s}_{N',m'}\big|^{4},$$
so that $A_4(N,N')=\int_{\mathcal{M}}G(\Psi')\cdot\Psi\,d\lambda$. 

\begin{lemma}[The Fourier series of the fourth power of the structured component]\label{cone-struct}
Let $\psi:=\psi_{N,1}$ be the kernel of Theorem \textup{\ref{FHdecom}(i)} and put
$H_m:=\big|\xi^{s}_{N,m}\big|^{2}$, so that $\Psi=\mathbb{E}_m H_m^2$. Then for every $m$,
$$H_m=\sum_{t\in\mathbb{Q}^{+}}h_m(t)\,\xi_t,\qquad h_m(t)\geq 0,\quad h_m(t)=h_m(1/t),$$
with finitely many non-zero terms. Consequently,
$$\Psi=\sum_{t\in\mathbb{Q}^{+}}d_t\,\xi_t,\qquad d_t\geq0,\quad d_t=d_{1/t},$$
and $\Psi$ is real valued and non-negative on $\mathcal{M}$. The same holds for $\Psi'$.
\end{lemma}

\begin{proof}
By Theorem \ref{FHdecom}(i) and the definition of the convolution on
$\mathbb{Z}/\tilde{N}\mathbb{Z}$,
$$\xi^{s}_{N,m}(\chi)=\chi_N\ast\psi(m)
=\underset{k\in\mathbb{Z}/\tilde{N}\mathbb{Z}}{\mathbb{E}}\psi(m-k)\chi_N(k)
=\sum_{k=1}^{N}a_m(k)\,\xi_k(\chi),\qquad a_m(k):=\frac{\psi(m-k)}{\tilde{N}}\geq0,$$
since $\psi$ is a kernel, hence non-negative. As every $\chi\in\mathcal{M}$ is completely
multiplicative and unimodular it extends to a homomorphism $\mathbb{Q}^{+}\to S^{1}$ by
$\chi(p/q)=\chi(p)\overline{\chi(q)}$, so that $\xi_p\overline{\xi_q}=\xi_{p/q}$. Hence
$$H_m=\xi^{s}_{N,m}\cdot\overline{\xi^{s}_{N,m}}
=\sum_{p,q=1}^{N}a_m(p)a_m(q)\,\xi_{p/q}
=\sum_{t}h_m(t)\xi_t,\qquad h_m(t)=\sum_{p/q=t}a_m(p)a_m(q)\geq0,$$
and $h_m(t)=h_m(1/t)$ upon interchanging $p$ and $q$. Squaring, the coefficient sequence of
$H_m^2$ is the multiplicative convolution of $h_m$ with itself, hence again non-negative and
symmetric, and averaging over $m$ preserves both properties. Finally $H_m=|\xi^s_{N,m}|^2\geq0$,
so $\Psi=\mathbb{E}_mH_m^2\geq0$. The proof for $\Psi'$ is symmetric.
\end{proof}

We let $\chi_{\mathrm{triv}}:\Q^+\rightarrow S^1$ denote the trivial character sending every $t\in \Q^+$ to $1$.
\begin{lemma}[The total mass is bounded below]\label{mass}
With $C:=\frac{1}{20^{4}l^{4}}$ we have
$$\Psi(\chi_{\mathrm{triv}})=\sum_{t\in\mathbb{Q}^{+}}d_t\ \geq\ C\ >\ 0,$$
and likewise $\Psi'(\chi_{\mathrm{triv}})\geq C$.
\end{lemma}

\begin{proof}
Since $\xi_t(\chi_{\mathrm{triv}})=\chi_{\mathrm{triv}}(t)=1$ for every $t\in\mathbb{Q}^{+}$, evaluating the
expansion of Lemma \ref{cone-struct} at the trivial character $\chi_{\mathrm{triv}}$ gives
$\Psi(\chi_{\mathrm{triv}})=\sum_t d_t$. On the other hand
$\Psi(\chi_{\mathrm{triv}})=\mathbb{E}_m\big|\xi^{s}_{N,m}(\chi_{\mathrm{triv}})\big|^{4}
=\mathbb{E}_m\big|{\chi_{\mathrm{triv},N,m}^s}\big|^{4}$, where
$\chi_{\mathrm{triv},N}^{s}=\mathbf \chi_{\mathrm{triv},N}\ast\psi_{N,1}$ and $\chi_{\mathrm{triv},N}=1_{[N]}$. As
$\psi_{N,1}$ is a kernel, $\chi_{\mathrm{triv},N,m}^{s}\geq0$ for every $m$ and
$$\underset{m\in\mathbb{Z}/\tilde{N}\mathbb{Z}}{\mathbb{E}}\chi_{\mathrm{triv},N,m}^{s}
=\underset{m}{\mathbb{E}}\ \underset{k}{\mathbb{E}}\,\psi_{N,1}(m-k)1_{[N]}(k)
=\underset{k}{\mathbb{E}}\,1_{[N]}(k)\underset{m}{\mathbb{E}}\,\psi_{N,1}(m-k)
=\frac{N}{\tilde{N}} ,$$
since $\psi_{N,1}$ has average $1$. By Jensen's inequality applied to the convex function
$t\mapsto t^{4}$ on $[0,\infty)$,
\begin{equation}\label{Jensen}
\Psi(\chi_{\mathrm{triv}})=\underset{m}{\mathbb{E}}\big|\chi_{\mathrm{triv},N,m}^{s}\big|^{4}
\ \geq\ \Big(\underset{m}{\mathbb{E}}\,\chi_{\mathrm{triv},N,m}^{s}\Big)^{4}
=\Big(\frac{N}{\tilde{N}}\Big)^{4}\ \geq\ \frac{1}{20^{4}l^{4}}=C,
\end{equation}
where we used $\tilde{N}\leq 20lN$. The argument for $\Psi'$ is identical.
\end{proof}

By the previous lemma we can set $T_{\Psi}:=\sum_{t\in \Q^{+}} d_t T_t$, and $M:=T_{\Psi}/\Psi(\chi_{\mathrm{triv}})$ and $M':=T_{\Psi'}/\Psi'(\chi_{\mathrm{triv}})$. We write $M'$ rather than $N$, because $N$ is already used in our computations.
\begin{lemma}\label{psd-markov}
In the setting above, $M$ and $M'$ are PSD doubly-stochastic operators.
\end{lemma}

\begin{proof}
The sums are finite by Lemma \ref{cone-struct}. If $g\geq0$ then $T_tg\geq0$ for every $t$
and $d_t\geq0$, so $T_\Psi g\geq0$. Since $T_t\mathbf 1_X=\mathbf 1_X$ we get
$T_\Psi\mathbf 1_X=\sum_t d_t T_t\bm{1}_X = (\sum_t d_t) \bm{1}_X = \Psi(\chi_{\mathrm{triv}})\cdot \mathbf 1_X$, and since each $T_t$ preserves $\mu$ we get
$\int T_\Psi f\,d\mu=\Psi(\chi_{\mathrm{triv}})\cdot\int f\,d\mu$. Dividing by $\Psi(\chi_{\mathrm{triv}})$ gives $M\mathbf 1_X=\mathbf 1_X$
and $\int Mf\, d\mu=\int f\, d\mu$. Each $T_t$ is unitary with $T_t^{*}=T_{t^{-1}}$, so by the symmetry
$d_t=d_{1/t}$ of Lemma \ref{cone-struct},
$$T_\Psi^{*}=\sum_t d_t T_{t^{-1}}=\sum_t d_{1/t}T_t=T_\Psi .$$
For positive semidefiniteness, write $T_c:=\sum_t c(t)T_t$ for a finitely supported $c$.
Since $h_m$ is real and symmetric and the coefficient sequence of $H_m^{2}$ is the
convolution of $h_m$ with itself,
$$T_{h_m}T_{h_m}^{*}
=\sum_{t_1,t_2}h_m(t_1)h_m(t_2)\,T_{t_1t_2^{-1}}
=T_{H_m^{2}} .$$
Averaging over $m$ gives $T_\Psi=\mathbb{E}_m\,T_{h_m}T_{h_m}^{*}$, whence for every
$g\in L^2(\mu)$
$$\langle T_\Psi g,g\rangle=\underset{m}{\mathbb{E}}\,\big\|T_{h_m}^{*}g\big\|_{2}^{2}\geq0 .$$
The argument for $M'$ is identical.
\end{proof}

With this notation, \eqref{spec2gen} applied with $f=g=h=1_A$, together with the
self-adjointness of $T_{\Psi'}$, gives
\begin{equation}\label{A4identity}
A_4(N,N')=\int_X f\cdot\big(T_\Psi f\big)\cdot\big(T_{\Psi'}f\big)\,d\mu
=\Psi(\chi_{\mathrm{triv}})\Psi'(\chi_{\mathrm{triv}})\int_X f\,(Mf)(M'f)\,d\mu .
\end{equation}

We are now in a position to apply Theorem \ref{conj} of Section \ref{positivity}.

\begin{lemma}\label{A4lower}
$A_4(N,N')\geq C^{2}\delta^{3}$, where $C=\frac{1}{20^4l^4}$.
\end{lemma}

\begin{proof}
By \eqref{A4identity}, Lemma \ref{psd-markov} and Theorem \ref{conj},
$$A_4(N,N')=\Psi(\chi_{\mathrm{triv}})\,\Psi'(\chi_{\mathrm{triv}})\int_X f\,(Mf)(M'f)\,d\mu
 \geq\ \Psi(\chi_{\mathrm{triv}})\,\Psi'(\chi_{\mathrm{triv}})\delta^{3}\geq\ C^{2}\delta^{3},$$
the last inequality by Lemma \ref{mass}.
\end{proof}

\subsection{Concluding the proof}
We can now choose the constants. Recall from \eqref{Gbound} that $\|G\|_{op}\leq 1$.
Combining \eqref{A2 - A3}, \eqref{A3-A4} and Lemma \ref{A4lower},
$$A_3(N,N')\geq A_4(N,N')-|A_3-A_4|\geq C^2\delta^{3}-16\varepsilon\|G\|_{op},\qquad
A_2(N,N')\geq C^2\delta^{3}-16\varepsilon\|G\|_{op}-8\sqrt{\varepsilon}\|G\|_{op}.$$
Recall $\varepsilon=c_3\delta^4$. By \eqref{Gbound},
$$16\varepsilon\|G\|_{op}\leq 16c_3\delta^{4},
\qquad 8\sqrt{\varepsilon}\|G\|_{op}\leq 8\sqrt{c_3}\,\delta^{2},$$
so choosing
$$c_3\ :=\ \min\Big\{1,\ \frac{C^{2}}{64\delta},\ \frac{C^{4}\delta^{2}}{1024}\Big\}$$
makes each of the two error terms at most $\tfrac14C^2\delta^{3}$, and hence
\begin{equation}\label{A2lower}
A_2(N,N')\ \geq\ \tfrac12\,C^{2}\delta^{3}.
\end{equation}
By \eqref{A1 geq A2} and \eqref{A2lower},
$$A_1(N,N')\ \geq\ \frac{\varepsilon^{2}}{500l^{2}\,Q^{2}R^{2}}\cdot\tfrac12C^{2}\delta^{3}
=\frac{c_3^{2}C^{2}\,\delta^{11}}{1000l^2\,Q^{2}R^{2}} .$$
Finally, by \eqref{A-A1} and \eqref{Gbound}, and since
$F(Q,R,\varepsilon)^{-1/2}=\frac{\sqrt{\varepsilon}}{c_4Q^{2}R^{2}}$,
$$|A(N,N')-A_1(N,N')|\ \leq\ \frac{8\,c_2\sqrt{c_3}\,\delta^{2}}{c_4\,Q^{2}R^{2}}
+\frac{16}{\tilde{N}}+\frac{16}{\tilde{N}'} .$$
Both main terms carry the same factor $Q^{-2}R^{-2}$, so $Q$ and $R$ cancel. Choosing
$$c_4\ :=\ \frac{20^4\,c_2\,l^{2}}{c_3^{3/2}\,C^{2}\,\delta^{9}},$$
the displayed error is at most half of the lower bound for $A_1$, and therefore
$$A(N,N')\ \geq\ \frac{c_3^{2}C^{2}\delta^{11}}{2000\,l^{2}\,Q^{2}R^{2}}
-\frac{16}{\tilde{N}}-\frac{16}{\tilde{N}'}\ \geq\ \frac{c_3^{2}C^{2}\delta^{11}}{2000\,l^{2}\,K^{4}}
-\frac{16}{\tilde{N}}-\frac{16}{\tilde{N}'}.$$
By Remark \ref{pigeonhole} the integers $Q$ and $R$ do not depend on $N,N'$, and they are bounded by $K=K(\varepsilon,F)$, hence by a constant depending only on $l$ and $\delta$. Since these occur in the denominator, the last bound is uniform in $N,N'$. Since $\tilde{N}\geq 10lN$ and $\tilde{N}'\geq 10lN'$, the two error terms tend to $0$, so
$$A(N,N') \geq \frac{c_3^{2}C^{2}\delta^{11}}{4000\,l^{2}\,K^{4}} > 0
\qquad\text{for all sufficiently large } N,N'\in\mathcal{N},$$
and $\mathcal{N}$ is infinite, so this holds for arbitrarily large $N,N'$.
By Lemma \ref{liminf'} this proves \eqref{liminf} with $\theta:=\frac{c_3^{2}C^{2}\delta^{11}}{4000\,l^{2}\,K^{4}}$, and with it Theorem~\ref{multrec} and Theorem~\ref{PR}. $\square$

\bibliographystyle{abbrv}
\bibliography{bibliography}
\end{document}